\documentclass[11pt,twoside]{article}
\usepackage{amsfonts,amsmath,amssymb}
\usepackage{mathrsfs}
\usepackage{picinpar}

\newcommand{\lbl}[1]{\label{#1}}

\newtheorem{theo}{Theorem}[section]
\newtheorem{prop}{Proposition}[section]
\newtheorem{lem}{Lemma}[section]
\newtheorem{remark}{Remark}[section]
\newtheorem{cor}{Corollary}[section]
\newtheorem{defi}{Definition}[section]
\newcommand{\be}{\begin{equation}}
\newcommand{\ee}{\end{equation}}
\newcommand\bes{\begin{eqnarray}} \newcommand\ees{\end{eqnarray}}
\newcommand{\bess}{\begin{eqnarray*}}
\newcommand{\eess}{\end{eqnarray*}}
\newcommand\bedd{\bes\left\{\begin{array}{ll}\medskip}
\newcommand\eedd{\end{array}\right.\ees}

\newcommand\ep{\varepsilon}
\newcommand\kk{\left}
\newcommand\rr{\right}
\newcommand\nm{\nonumber}
\newcommand\dd{\displaystyle}
\newcommand\vp{\varphi}

\newcommand\mR{\mathbb{R}}
\newcommand\ud{\underline}
\newcommand\od{\overline}
\newcommand\lm{\lambda}

\newcommand\yy{\infty}

\newcommand\qqq{\eqref}

\newcommand\ff{, \ \ \forall \ }
\begin{document}
  \pagestyle{myheadings}
\thispagestyle{empty}

\begin{center}{\Large Free boundary problems of the competition model with}\\[2mm]
{\Large sign-changing coefficients in heterogeneous time-periodic
  environment}\footnote{This work was supported by NSFC Grant 11371113}\\[5mm]
 {\Large  Mingxin Wang\footnote{E-mail: mxwang@hit.edu.cn; Tel: 86-15145101503; Fax:
 86-451-86402528}}\\[2mm]
{Natural Science Research Center, Harbin Institute of Technology, Harbin 150080, PR
China}
\end{center}

\begin{quote}
\noindent{\bf Abstract.} In this paper we consider two kinds of free boundary problems for the diffusive competition model in the heterogeneous time-periodic environment, in which the variable intrinsic growth rates of invasive and native species may change signs and be ``very negative" in a ``suitable large region" (see the conditions {\bf(A)}
and {\bf(H1)}). The main purpose is to
understand the dynamical behavior of the two competing species spreading
via a free boundary. We study the spreading-vanishing dichotomy, long time behavior of solution, sharp criteria for spreading and vanishing, and estimates of the asymptotic spreading speed of the free boundary. Moreover, the existence of positive solutions to a $T$-periodic boundary value problem in half line, associated with our free boundary problems, is obtained.


 \end{quote}

 \section{Introduction}
 \setcounter{equation}{0} {\setlength\arraycolsep{2pt}

In the natural world, the following phenomenon will happen constantly:

$\bullet$\, There is one kind of native species in an area (initial habitat). At some time (initial time), a new or invasive species (competitor) enters this district.

A typical mathematical model describing the interaction (competition)
between invasive and native species is the following competitive model
 \bes\left\{\begin{array}{lll}
 u_t-d_1 u_{xx}=u\big(a(t,x)-c(t,x)u-k(t,x)v\big),\\[1mm]
 v_t-d_2 v_{xx}=v\big(b(t,x)-m(t,x)v-h(t,x)u\big),
 \end{array}\right.\lbl{1.1}\ees
where $u(t,x)$ and $v(t,x)$ represent the population densities of the invasive and native species, respectively; $d_1,\,d_2>0$ are their diffusion (dispersal) rates;
$a(t,x)$, $b(t,x)$ denote their intrinsic growth rates; $c(t,x)$, $m(t,x)$ are the
intraspecific and $k(t,x)$, $h(t,x)$ the interspecific competition rates. The system
(\ref{1.1}), as a model describing the spreading, persistence and extinction of two
competing species in the heterogeneous environment, has received an astonishing
amount of attention, please refer to \cite{CC2, CLL, DHMP, HMP, LL, Ni} for example.
When the functions $a,b,c,m,k$ and $h$ are positive constants, to describe the invasion and spreading phenomenon, there have been many interesting studies on positive traveling
waves and asymptotic spreading speed of (\ref{1.1}), see, for example \cite{GL, LLW,
Zhaox} and the references cited therein.

In the natural world, for most animals and plants, their birth and death rates will
change with seasons, so the intrinsic growth rates $a(t,x),b(t,x)$ and then the
intraspecific competition rates $c(t,x)$, $m(t,x)$ and interspecific competition
rates $k(t,x)$, $h(t,x)$ should be time-periodic functions. Especially, in the
winter of severe cold and cold zones, animals cannot capture enough food to feed
upon and will not breed, seeds cannot germinate and buds cannot grow above ground,
so their birth rates are zero. In the meantime, their death rate will be greater. Therefore, in some periods and some areas, the intrinsic growth rates $a(t,x)$ and
$b(t,x)$ may be negative.

In general, the invasive and native species will have a tendency to emigrate from
the boundary to obtain their new habitats, i.e., they will move outward along the
unknown curve (free boundary) as time increases. In order to simplify the
mathematics, in this paper we only consider the one dimensional case and assume that
the left boundary is fixed: $x=0$. Moreover, we take $c(t,x),\,m(t,x),\,k(t,x)$ and
$h(t,x)$ are positive constants, and by the suitable rescaling we may think that
$c(t,x)=m(t,x)=1$. It should be emphasized that for the higher dimensional and radially symmetric case, when $c,\,m,\,k$ and $h$ are functions of $(t,x)$
and have positive lower and upper bounds, the methods used in this
paper are still valid and related results remains hold.

If, at the initial time, the range occupied by native species is not very large and the
invasive species has a wide distribution in this area, we can use the following free
boundary problem to model the above phenomenon:
  \bes\left\{\begin{array}{lll}
 u_t-d_1u_{xx}=u\big(a(t,x)-u-kv\big), &t>0, \ \ 0<x<s(t),\\[1mm]
 v_t-d_2v_{xx}=v\big(b(t,x)-v-hu\big),&t>0, \ \ 0<x<s(t),\\[1mm]
 B_1[u]=B_2[v]=0,\ \ &t>0, \ \ x=0,\\[1mm]
 u=v=0, \ s'(t)=-\mu(u_x+\rho v_x),\ \ &t>0, \ \ x=s(t), \\[1mm]
 u(0,x)=u_0(x), \ \ v(0,x)=v_0(x),& 0\leq x\leq s_0,\\[1mm]
 s(0)=s_0.
 \end{array}\right.\lbl{1.2}\ees

If, at the initial time, the distribution range of native species is very large (it can be considered to be half line) and the invasive species has a local distribution in
this area, we shall use the following free boundary problem to describe the above
phenomenon:
 \bes\left\{\begin{array}{lll}
 u_t-d_1u_{xx}=u\big(a(t,x)-u-kv\big), &t>0, \ \ 0<x<s(t),\\[1mm]
 u(t,x)\equiv 0,&t>0, \ \ x\ge s(t),\\[1mm]
 v_t-d_2v_{xx}=v\big(b(t,x)-v-hu\big),&t>0, \ \ 0<x<\infty,\\[1mm]
 B_1[u]=B_2[v]=0,\ \ &t>0, \ \ x=0,\\[1mm]
 u=0, \ s'(t)=-\mu u_x,\ \ &t>0, \ \ x=s(t), \\[1mm]
 u(0,x)=u_0(x), \ 0\leq x\leq s_0; \ \ v(0,x)=v_0(x),&x\ge 0,\\[1mm]
 s(0)=s_0.
 \end{array}\right.\lbl{1.3}\ees

In the above two problems, $B_1[u]=\alpha_1u-\beta_1u_x$,
$B_2[v]=\alpha_2v-\beta_2v_x$, $\alpha_i,\beta_i$ are nonnegative constants and
satisfy $\alpha_i+\beta_i=1$; $x=s(t)$ is the free boundary to be determined;
positive constant $s_0$ is the initial boundary or survival range; positive constants
$\mu$ and $\mu\rho$, the expansion capacities, are the ratios of the expansion speed of the free boundary relative to population gradients at the expanding front, those  describe the abilities of species to transmit and dispersal in the new habitat and can also be considered as the ``moving parameters".

Throughout this paper, we assume that
\vspace{-1mm}\begin{quote}
 {\bf(H)}\, Functions $a,\,b\in\big(C^{\frac{\nu}2,\nu}\cap
 L^\infty\big)([0,\infty)\times[0,\infty))$ for some $\nu\in(0,1)$,
 and are $T$-periodic in time $t$ for some $T>0$ and positive somewhere in $[0,T]\times[0,\infty)$;
  \end{quote}\vspace{-1mm}
and the initial functions $u_0(x),v_0(x)$ satisfy

$\bullet$\, $u_0,\,v_0\in C^2([0,s_0])$, $u_0,\,v_0>0$ in $(0,s_0)$,
$B_1[u_0](0)=u_0(s_0)=0$ and $B_2[v_0](0)=v_0(s_0)=0$ for the problem (\ref{1.2});

$\bullet$\, $u_0\in C^2([0,s_0])$, $v_0\in C^2([0,\infty))$, $u_0>0$ in
$(0,s_0)$, $v_0>0$ in $(0,\infty)$, $B_1[u_0](0)=u_0(s_0)=0$ and $B_2[v_0](0)=0$ for
the problem (\ref{1.3}).

\vskip 2pt Some problems associated with (\ref{1.2}) and (\ref{1.3}) have been studied recently. When the functions $a,\,b$ are positive constants, the problem (\ref{1.2}) has been studied by Guo \& Wu \cite{GW} and Wang \& Zhao \cite{WZjdde} for the case that $\alpha_1=\alpha_2=0$, and by Wang \& Zhao \cite{WZjdde} for the case that
$\beta_1=\beta_2=0$. When the functions $a,\,b$ are positive constants, or $a,\,b\in \big(C^1\cap L^\infty\big)([0,\infty))$ are independent of time $t$ and strictly positive, the problem (\ref{1.3}) with $\alpha_1=\alpha_2=0$ has been
studied by Du \& Lin \cite{DL2} and Wang \& Zhang \cite{WZ15} for the higher dimension and radially symmetric
case. The complete conclusions
about spreading-vanishing dichotomy, sharp criteria for spreading and vanishing,
long time behaviour of $(u,v)$ and asymptotic spreading speed of the free boundary have been obtained in \cite{GW, WZjdde, DL2}.
Some  free boundary problems of diffusive prey-predator model with positive constant
coefficients has been studied by Wang, Zhang \& Zhao \cite{Wjde14, Wcnsns15, WZhang, WZ, ZW}.

In the absence of a native species, namely $v\equiv 0$, both problems (\ref{1.2})
and (\ref{1.3}) reduce to the following diffusive logistic problem with a free boundary
\bes
 \left\{\begin{array}{lll}
 u_t-d_1u_{xx}=u\big(a(t,x)-u\big), \ \ &t>0,\ \ 0<x<s(t),\\[1mm]
 B_1[u](t,0)=0,\ \ u(t,s(t))=0,\ \ &t>0,\\[1mm]
 s'(t)=-\mu u_x(t,s(t)),&t>0,\\[1mm]
 s(0)=s_0,\ \ u(0,x)=u_0(x),&0\le x\le s_0,
 \end{array}\right.\label{1.5}
 \ees
which may be used to describe the spreading of a new or invasive species and has been studied by the author in \cite{Wperiodic15} recently. When the function
$a$ has positive lower and upper bounds, i.e., there exist positive constants
$\kappa_1$, $\kappa_2$ such that $\kappa_1\le a(t,x)\le \kappa_2$, Du, Guo \& Peng in \cite{DGP} have studied the problem (\ref{1.5}) with $\alpha_1=0$
for the higher dimension and radially symmetric case. The spreading-vanishing
dichotomy, sharp criteria for spreading and vanishing and asymptotic spreading speed
of the free boundary have been obtained in \cite{DGP, Wperiodic15}.

When $a=a(x)$ is independent of the time $t$ and changes sign, the problem
(\ref{1.5}) was studied by Zhou \& Xiao \cite{ZX} and Wang \cite{Wjde15}. When $a=a(x)$ has positive
lower and upper bounds, some similar problems to (\ref{1.5}) has been studied
systematically.  When $a$ is a positive constant, the problem (\ref{1.5}) was
investigated earlier by Du \& Lin \cite{DLin} for $\alpha_1=0$ and by Kaneko \&
Yamada \cite{KY} for $\beta_1=0$. Du, Guo \& Liang \cite{DG, DLiang} discussed the
higher dimensional and radially symmetric case ($\alpha_1=0$). The non-radial case in higher dimensions was treated by Du \& Guo \cite{DG1}. Peng \&
Zhao \cite{PZ} studied a free boundary problem of the diffusive logistic model with
seasonal succession. They considered that the species does not migrate and stays in
a hibernating status in bad season. The evolution of the species obeys Malthusian's
equation $u_t=-\delta u$ in  bad season, and obeys the diffusive logistic equation
with positive constant coefficients in good season. Instead of $u(a-bu)$ by a general function $f(u)$, Du, Matsuzawa \& Zhou \cite{DMZ}, Kaneko \cite{Ka} and Du \& Lou \cite{DLou} investigated the corresponding free boundary problems.

The main aim of this paper is to study the dynamical properties of (\ref{1.2}) and
(\ref{1.3}). We first state the global existence, uniqueness, regularity and
estimate of solution $(u,v,s)$.
  \begin{theo}\lbl{th2.1}\, The problem {\rm(\ref{1.2})} has a unique global solution
$(u,v,s)$ in time and satisfies
 \bes
  u,v\in C^{1+\frac{\nu} 2,2+\nu}(D_\infty),\ \ s\in
  C^{1+\frac{1+\nu}2}((0,\infty)),
  \lbl{2.1}\ees
where $D_\infty=\big\{t>0,\, x\in[0,s(t)]\big\}$. Furthermore, there exist constants
$M=M\kk(\|a,b,u_0,v_0\|_\infty\rr)>0$ and
$C=C\kk(\mu,\,\|a,b,u_0,v_0\|_\infty\rr)>0$, such that
  \bes
  0<u(t,x),\ v(t,x)\leq M, \ \ 0<s'(t)\leq \mu M, \ \ \forall \  t> 0,\
  0<x<s(t)\lbl{2.2}\ees
and
 \bes\left\{\begin{array}{ll}
  \|u(t,\cdot),\,v(t,\cdot)\|_{C^1[0,\,s(t)]}\leq C, \ \ &\forall \ t\ge
  1, \\[1.5mm]
  \|s'\|_{C^{{\nu}/2}([n+1,n+3])}\leq C, \ \ &\forall \ n\geq 0.
  \end{array}\right.\lbl{2.3}\ees
\end{theo}
 \begin{theo}\lbl{th6.1}\, The problem {\rm(\ref{1.3})} has a unique global solution
$(u,v,s)$ in time and
 \bess
  u\in C^{1+\frac{\nu} 2,2+\nu}(D_\infty),\ \ v\in C^{1+\frac{\nu}
  2,2+\nu}((0,\infty)\times[0,\infty)), \ \ s\in
  C^{1+\frac{1+\nu}2}((0,\infty)).
  \eess
Moreover, $(u,v,s)$ satisfies the estimates
  \bess
  0<u(t,x),\ v(t,x)\leq M, \ \ 0<s'(t)\leq \mu M, \ \ \forall \  t,\,x>
  0\eess
and
 \bes\left\{\begin{array}{ll}
  \|u(t,\cdot)\|_{C^1([0,\,s(t)])}, \ \|v(t,\cdot)\|_{C^1([0,\infty))}\leq
  C, \ \ &\forall \  t\ge 1,\\[1.5mm]
  \|s'\|_{C^{\nu/2}([n+1,n+3])}\leq C, &\forall \  n\geq 0,
  \end{array}\right.\lbl{1.8}\ees
where $D_\infty$, $M$ and $C$ are same as those of Theorem $\ref{th2.1}$.
\end{theo}

Proofs of Theorems \ref{th2.1} and \ref{th6.1} are essentially parallel to that of \cite{DL2, GW, Wperiodic15}. For the global existence and uniqueness of $(u,v,s)$, please refer to \cite[Theorem 2.1]{DL2} and \cite[Theorem 1]{GW}; for the regularities and estimates of $(u,v,s)$, please refer to \cite[Theorem 2.1]{Wperiodic15}. The details are omitted here. We remark that the uniform estimates (\ref{2.3}) and (\ref{1.8}) allow us to deduce that $s'(t)\to 0$ when $s_\infty<\infty$ and play a key role for determining the vanishing phenomenon.

\vskip 2pt It follows from Theorems \ref{th2.1} and \ref{th6.1} that $s(t)$ is monotonically increasing. There exists $s_\infty\in(0,\infty]$ such that $\dd\lim_{t\to\infty} s(t)=s_\infty$.

In order to study the long time behavior of solution
and spreading phenomenon, in Section 2 we investigate a
stationary problem: the $T$-periodic boundary value problem in half line associated with the free boundary problems  (\ref{1.2}) and (\ref{1.3}). As most properties of (\ref{1.2}) and (\ref{1.3}) are similar, we first deal with the problem (\ref{1.2}) meticulously in Sections 3 and 4, and briefly discuss the problem
(\ref{1.3}) in Section \ref{s5}. In Section 3, we shall derive the spreading-vanishing dichotomy of (\ref{1.2}):

Either
\vspace{-2mm}\begin{itemize}
\item[(i)]spreading: $s_\yy=\infty$ and
\bess
 & U_*(t,x)\le\dd\liminf_{n\to\infty}u(t+nT,x), \ \
 \limsup_{n\to\infty}u(t+nT,x)\leq U^*(t,x),&\\
 &V_*(t,x)\le\dd\liminf_{n\to\infty}v(t+nT,x), \ \
 \limsup_{n\to\infty}v(t+nT,x)\leq V^*(t,x),&
 \eess
uniformly in $[0, T]\times[0,L]$ for any $L>0$, where $(U^*,V_*)$ and
$(U_*, V^*)$ are positive $T$-periodic solutions of (\ref{3.1}) which will
be given in Theorem $\ref{th3.1}$;\vspace{-2mm}\end{itemize}
or
\vspace{-2mm}\begin{itemize}
\item[(ii)] vanishing: $s_\infty<\infty$ and
 $\dd\lim_{t\to\infty}\|u(t,\cdot),\,v(t,\cdot)\|_{C([0,s(t)])}=0$.
 \vspace{-2mm}\end{itemize}
In Section 4, the criteria for spreading and vanishing of the problem
$(\ref{1.2})$ will be  established.

\section{Positive solutions of a $T$-periodic boundary value problem in half line}
 \setcounter{equation}{0} {\setlength\arraycolsep{3pt}

To discuss the long time behavior of solution and spreading phenomenon, we should
study the following $T$-periodic boundary value problem in half line
 \bes\left\{\begin{array}{ll}
 U_t-d_1U_{xx}=U\big(a(t,x)-U-k V\big), \ \ &0<t<T, \
 0<x<\yy,\\[1mm]
 V_t-d_2V_{xx}=V\big(b(t,x)-V-h U\big), \ \ &0<t<T, \
 0<x<\yy,\\[1mm]
 B_1[U](t,0)=B_2[V](t,0)=0,  &0\le t\le T,\\[1mm]
 U(0,x)=U(T,x),\ \ V(0,x)=V(T,x), &0\le x<\yy.
 \end{array}\right. \lbl{3.1}\ees

\subsection{Preliminaries}

In order to facilitate applications, in this subsection we state some known
results: the comparison principle, properties of the principal eigenvalue of the
$T$-periodic eigenvalue problem and conclusions of the diffusive logistic equation in the heterogeneous time-periodic environment.
 \begin{lem}\lbl{l3.1}{\rm(}Comparison principle{\rm)}\, Let $\tau, \ell>0$. Assume
that
 \[\od w,\,\ud w,\,\od z,\,\ud z\in C([0,\tau]\times[0,\ell])\cap
 C^{0,1}([0,\tau]\times[0,\ell))\cap C^{1,2}((0,\tau]\times(0,\ell)),\]
and are nonnegative functions. If $(\od w,\ud z)$ and $(\ud w,\od z)$ satisfy
 \bess\left\{\begin{array}{ll}
 \od w_t-d_1\od w_{xx}\ge\od w\big(a(t,x)-\od w-k\ud z\big), \ \ &0<t\le\tau, \
 0<x<\ell,\\[1mm]
 \ud z_t-d_2\ud z_{xx}\le\ud z\big(b(t,x)-\ud z-h\od w\big), \ \ &0<t\le\tau, \
 0<x<\ell,
 \\[1mm]
   \ud w_t-d_1\ud w_{xx}\le\ud w\big(a(t,x)-\ud w-k\od z\big), \ \ &0<t\le\tau, \
   0<x<\ell,\\[1mm]
 \od z_t-d_2\od z_{xx}\ge\od z\big(b(t,x)-\od z-h\ud w\big), \ \ &0<t\le\tau, \
 0<x<\ell
 \end{array}\right.\eess
 and
 \bess\left\{\begin{array}{ll}
  B_1[\od w](t,0)\ge B_2[\ud z](t,0), &0\le t\le\tau,\\[1mm]
  B_1[\ud w](t,0)\le B_2[\od z](t,0), &0\le t\le\tau,\\[1mm]
 \od w(t,\ell)\ge\ud w(t,\ell), \ \ \ud z(t,\ell)\le\od z(t,\ell), \ \ &0\le
 t\le\tau,\\[1mm]
 \od w(0,x)\ge\ud w(0,x),\ \ \ud z(0,x)\le \od z(0,x), &0\le x\le \ell.
 \end{array}\right.\eess
Then we have
 \[\ud w\le\od w, \ \ \ud z\le\od z \ \ \ \mbox{in} \ \ [0,\tau]\times[0,\ell].\]
\end{lem}

{\bf Proof}. Note that the function pair $\big(u(a(t,x)-u-kv), v(b(t,x)-v-hu)\big)$
of $(u,v)$ is {\it quasimonotone nonincreasing} in $u,v\ge 0$. The desired result
can be deduced by the comparison principle for parabolic systems (see \cite{PW} or
\cite[section 4.2.2]{YLWW}). We omit the details here.\ \ \ \fbox{}

\vskip 2pt In the following we assume that $c,\,q\in\big(C^{\frac{\nu}2,\nu}\cap
L^\infty\big)([0,\infty)\times[0,\infty))$ and are $T$-periodic
functions in time $t$. Moreover, there exist positive constants $\ud q,\,\od q$ such that
$\ud q\leq q(t,x)\le\od q$ for all $t,x\ge 0$. Define $B[u]=\alpha u-\beta u_x$,
where $\alpha$ and $\beta$ are non-negative constants and satisfy $\alpha+\beta=1$.
 \begin{lem}{\rm(\cite[Lemma 3.2]{Wperiodic15})}\lbl{l3.3}\, Let $d,\ell>0$, $\od z,\,\ud
z \in C^{1, 2}([0,T]\times(0,\ell))\cap C^{0, 1}([0,T]\times[0,\ell])$. If $(\od
z,\ud z)$ satisfies
 \bess\left\{\begin{array}{ll}
  \od z_t-d\od z_{xx}\geq c(t,x)\od z-q(t,x)\od z^2,\ \ &0\le t\le T,\
  0<x<\ell,\\[1mm]
   \ud z_t-d\ud z_{xx}\le c(t,x)\ud z-q(t,x)\ud z^2,\ \ &0\le t\le T, \
   0<x<\ell,\\[1mm]
 B[\od z](t,0)\geq B[\ud z](t,0),\ \od z(t,\ell)\ge \ud z(t,\ell), \ \ &0\le t\le
 T,\\[1mm]
 \od z(0,x)=\od z(T,x),\ \ \ud z(0,x)=\ud z(T,x), &0\le x\le \ell.
 \end{array}\right.\eess
Then $\od z\ge\ud z$ in $[0,T]\times[0,\ell]$.
\end{lem}

For any given $d,\ell>0$, let $\lm_1(\ell;d,c)$ be the principal eigenvalue of the following $T$-periodic eigenvalue problem
 \bes\left\{\begin{array}{ll}
 \phi_t-d\phi_{xx}-c(t,x)\phi=\lm\phi, \ \ &0\le t\le T, \ 0<x< \ell,\\[1mm]
 B[\phi](t,0)=0, \ \ \phi(t, \ell)=0,\ \ &0\le t\le T, \\[1mm]
 \phi(0,x)=\phi(T,x),&0\le x\le \ell.
 \end{array}\right.\lbl{3.2}\ees
 \begin{prop}\lbl{p3.1}{\rm(\cite[Proposition 3.1]{Wperiodic15})}\, The principal
eigenvalue $\lm_1(\ell;d,c)$ is continuous with respect to $\ell,d$ and $c$, and
strictly decreasing in $c$ and $\ell$. Moreover, $\dd\lim_{\ell\to
0^+}\lm_1(\ell;d,c)=\infty$ and $\dd\lim_{d\to\infty}\lm_1(\ell;d,c)=\infty$.
 \end{prop}
  \begin{prop}\lbl{p3.2}{\rm(\cite[Proposition 3.2]{Wperiodic15})}\, Assume that the
function $c(t,x)$ satisfies
 \vspace{-1mm}\begin{quote}
 {\bf(A)}\, There exist $\varsigma>0$, $-2<r\le 0$, $m>1$ and $x_n$ satisfying
 $x_n\to\infty$ as $n\to\infty$, such that $c(t,x)\ge\varsigma x^r$ in
 $[0,T]\times[x_n,m x_n]$.
 \vspace{-1mm}\end{quote}
Then for any given $d>0$, there exists a unique $\ell_0=\ell_0(d)>0$ such that
$\lm_1(\ell_0;d,c)=0$. Hence, $\lm_1(\ell;d,c)<0$ for all $\ell>\ell_0$.
\end{prop}

Consider the following $T$-periodic boundary value problem of the diffusive logistic equation
in a bounded interval $(0,\ell)$:
\bes\left\{\begin{array}{ll}
 W_t-dW_{xx}=W\big(c(t,x)-q(t,x)W\big), \ \ &0\le t\le T, \
 0<x<\ell,\\[1mm]
 B[W](t,0)=0, \ \ W(t,\ell)=K, &0\le t\le T,\\[1mm]
 W(0,x)=W(T,x), &0\le x\le \ell.
  \end{array}\right. \lbl{3.4a}\ees
\begin{lem}\lbl{l3.2}{\rm(\cite[Lemma 3.3]{Wperiodic15})}\, Assume that $c(t,x)$ satisfies the condition {\bf(A)}. Then, for any given $\ell>\ell_0$ and $K\ge \|c\|_\yy/\ud q$, the
problem $(\ref{3.4a})$ has a unique positive solution.
\end{lem}

Now, let us consider the following initial-boundary value problem and $T$-periodic
boundary value problem of the diffusive logistic equation in the half line:
 \bes\left\{\begin{array}{ll}
 w_t-d w_{xx}=w\big(c(t,x)-w\big), \ \ &t>0, \ 0<x<\yy,\\[1mm]
 B[w](t,0)=0,  &t> 0,\\[1mm]
 w(0,x)=w_0(x), &0\le x<\yy
  \end{array}\right. \lbl{3.3}\ees
and
 \bes\left\{\begin{array}{ll}
 W_t-dW_{xx}=W\big(c(t,x)-W\big), \ \ &0\le t\le T, \
 0<x<\yy,\\[1mm]
 B[W](t,0)=0,  &0\le t\le T,\\[1mm]
 W(0,x)=W(T,x), &0\le x<\yy,
  \end{array}\right. \lbl{3.4}\ees
where $w_0(x)$ is a bounded nontrivial and nonnegative continuous function. For the convenience to write, we first a definition.
 \begin{defi}\lbl{d3.1}\, Let $r$ be a constant and satisfy $-2<r\le 0$,
$c\in\big(C^{\frac{\nu}2,\nu}\cap
L^\infty\big)([0,\infty)\times[0,\infty))$ be a $T$-periodic function
in time $t$. We call that $c$ belongs to the class ${\cal C}_r(T)$ if there exist two
positive $T$-periodic functions $c_\infty(t),c^\infty(t)\in C^{\nu/2}([0,T])$, such
that
 \bess
 c_\infty(t)\le\liminf_{x\to\infty}\frac{c(t,x)}{x^r}, \ \ \
 \limsup_{x\to\infty}\frac{c(t,x)}{x^r}\le c^\infty(t) \ \ \ \mbox{uniformly in } \
 [0,T].
 \eess
\end{defi}

It is easy to see that if $c\in {\cal C}_r(T)$ for some $-2<r\le 0$, then $c$
satisfies the condition {\bf(A)}.
 \begin{prop}\lbl{p3.3}\,Assume that $c\in {\cal C}_r(T)$. Then the $T$-periodic
boundary value problem $(\ref{3.4})$ has a unique positive solution $W\in
C^{1+\frac{\nu}2, 2+\nu}([0,T]\times[0,\infty))\cap{\cal C}_r(T)$, and satisfies
 \bes
 \min_{[0,T]}c_\infty(t)\le\liminf_{x\to\infty}\frac{W(t,x)}{x^r},\ \ \
 \limsup_{x\to\infty}\frac{W(t,x)}{x^r}\le \max_{[0,T]}c^\infty(t)
  \lbl{3.5}\ees
uniformly in $[0,T]$. Moreover, the solution $w$ of {\rm(\ref{3.3})} satisfies
 \bes
 \lim_{n\to\infty}w(t+nT)=W(t,x) \ \ \ \mbox{locally uniformly in } \
 [0,T]\times[0,\infty).\lbl{3.6}\ees

Especially, if $r=0$ then
 \bes
  w_\infty(t)\le\liminf_{x\to\infty}W(t,x),\ \ \
\limsup_{x\to\infty}W(t,x)\le w^\infty(t)
  \lbl{3.7}\ees
uniformly in $[0,T]$, where $w_\infty(t)$ and $w^\infty(t)$ are, respectively, the
unique positive solutions of the following $T$-periodic ordinary differential problems:
 \[w'(t)=w\big(c_\infty(t)-w\big), \ \ \ w(0)=w(T),\]
and
 \[w'(t)=w\big(c^\infty(t)-w\big), \ \ \ w(0)=w(T).\]
\end{prop}

{\bf Proof}. Existence and uniqueness of $W$ together with the estimate (\ref{3.5}) is just Theorem 4.2 of \cite{Wperiodic15}}. The limit (\ref{3.6}) is given by Theorem 4.3 of \cite{Wperiodic15}}. The estimate (\ref{3.7}) can be proved by the similar method to that of \cite[Theorem 1.4]{PD}. We omit the details here. \ \ \ \fbox{}
 \begin{prop}{\rm(}Comparison principle{\rm)}\lbl{p3.4}\, Assume that $-2<r_i\le 0$
and $c_i\in{\cal C}_{r_i}(T)$, $i=1,2$. Let $W_i\in C^{1+\frac{\nu}2,
2+\nu}([0,T]\times[0,\infty))$ be the unique positive solution of
$(\ref{3.4})$ with $c=c_i$. If $c_1\le c_2$, then $W_1\le W_2$ in
$[0,T]\times[0,\infty)$.
\end{prop}

{\bf Proof}.\, The existence and uniqueness of $W_i$ is guaranteed by
Proposition \ref{p3.3}. Since $c_i\in {\cal C}_{r_i}(T)$, we have that $c_i$ satisfies
the condition {\bf(A)}. In view of Proposition \ref{p3.2}, there exists
$\ell_0=\ell_0(d)\gg 1$ such that $\lm_1(\ell;d,c_i)<0$ for $\ell>\ell_0$ and
$i=1,2$.

For any given $\ell>\ell_0$. Utilizing Theorem 28.1 of \cite{Hess}, the
problem
  \bess\left\{\begin{array}{ll}
 W_t-d W_{xx}=W\big(c_i(t,x)-W\big), \ \ &0\le t\le T, \
 0<x<\ell,\\[1mm]
  B[W](t,0)=0, \ \ W(\ell, t)=0,  &0\le t\le T,\\[1mm]
 W(0,x)=W(T,x), &0\le x\le \ell
 \end{array}\right. \eess
has a unique positive $T$-periodic solution $W_{i\ell}(t,x)$. Since $c_1\le c_2$,
by Lemma \ref{l3.3},
 \bes W_{1\ell}(t,x)\le W_{2\ell}(t,x) \ \ \mbox{in} \ \
 [0,T]\times[0,\ell].
 \lbl{3.8}\ees

Obviously, $W_{i\ell}\le\|c\|_\infty$ by the maximum principle, and $W_{i\ell}$
is increasing in $\ell$ by Lemma \ref{l3.3}. Remember $W_i(t,x)$ is the unique
positive solution of (\ref{3.4}) with $c=c_i$. Make use of the regularity theory for
parabolic equations and compact argument, it can be proved that, for any given $L>0$, $W_{i\ell}\to W_i$ in $C^{1, 2}([0,T]\times[0,L])$ as $\ell\to\infty$, $i=1,2$. These facts combined with (\ref{3.8}) allow us to derive $W_1\le W_2$. The proof is
complete. \ \ \ \fbox{}

\subsection{Existence and properties of positive solutions to (\ref{3.1})}

We first state a condition.
 \vspace{-1mm}\begin{quote}\hspace{-2mm}
{\bf(H1)}\, The functions $a,b\in{\cal C}_r(T)$ for some $-2<r\le 0$. That is, there
exist positive $T$-periodic functions
$a_\infty(t),b_\infty(t),a^\infty(t),b^\infty(t)\in C^{\nu/2}([0,T])$, such that
 \bes\left\{\begin{array}{ll}
 a_\infty(t)\le\dd\liminf_{x\to\infty}\frac{a(t,x)}{x^r}, \ \ \
 \limsup_{x\to\infty}\frac{a(t,x)}{x^r}\le a^\infty(t),\\[3mm]
 b_\infty(t)\le\dd\liminf_{x\to\infty}\frac{b(t,x)}{x^r}, \ \ \
 \limsup_{x\to\infty}\frac{b(t,x)}{x^r}\le b^\infty(t)
 \end{array}\right.\lbl{1.4}\ees
uniformly in $[0,T]$.
  \vspace{-1mm}\end{quote}

When the condition {\bf(H1)} holds, we define
 \[\ud a_\infty=\min_{[0,T]}a_\infty(t), \ \ \bar a^\infty=\max_{[0,T]}a^\infty(t),
 \ \
 \ud b_\infty=\min_{[0,T]} b_\infty(t), \ \ \bar b^\infty=\max_{[0,T]}b^\infty(t).\]
  \begin{theo}\lbl{th3.1}\, Under the condition {\bf(H1)}, we assume further that
 \bes
 \ud b_\infty>h\od a^\infty, \ \ \ \ud a_\infty>k\od b^\infty.
 \lbl{3.11}\ees
Then there exist four positive $T$-periodic functions
$U^*,\,U_*,\,V^*,\,V_*\in C^{1+\frac{\nu}2,
2+\nu}([0,T]\times[0,\infty))$, such that both $(U^*,V_*)$ and
$(U_*,V^*)$ are positive solutions of $(\ref{3.1})$.  Moreover, any
positive solution $(U,V)$ of {\rm(\ref{3.1})} satisfies
\bes
 U_*\leq U\leq U^*, \  \ V_*\leq V\leq V^* \ \ \
 \mbox{in} \ \ [0,T]\times[0,\infty).
 \lbl{3.12}\ees

Especially, if $r=0$ in ${\bf(H1)}$, then any positive solution $(U,V)$ of
{\rm(\ref{3.1})} satisfies
 \bes\kk\{\begin{array}{ll}
 w_1(t)\le\dd\liminf_{x\to\infty}U(t,x),\ \ \
\limsup_{x\to\infty} U(t,x)\le w_2(t),\\[1mm]
z_1(t)\le\dd\liminf_{x\to\infty}V(t,x),\ \ \
\limsup_{x\to\infty} V(t,x)\le z_2(t)
 \end{array}\rr.\lbl{3.13}\ees
uniformly in $[0,T]$, where $w_2(t)$, $z_1(t)$, $ z_2(t)$ and $w_1(t)$ are the unique positive solutions of the following $T$-periodic ordinary
differential problems
 \bes
 &w_2'(t)=w_2\big(a^\infty(t)-w_2\big), \ \ \ w_2(0)=w_2(T),\nm&\\[1mm]
 &z_1'(t)=z_1\big(b_\infty(t)-h w_2(t)-z_1\big), \ \ \
 z_1(0)=z_1(T),\lbl{3.14}&\\[1mm]
 & z_2'(t)= z_2\big(b^\infty(t)- z_2\big), \ \ \  z_2(0)= z_2(T),\lbl{3.15}&
 \ees
and
  \[w_1'(t)=w_1\big(a_\infty(t)-k z_2(t)-w_1\big), \ \ \ w_1(0)=w_1(T),\]
respectively.
 \end{theo}
 \begin{remark}
 The condition $(\ref{3.11})$ is equivalent to $h<\ud b_\infty/\od a^\infty$ and $k<\ud a_\infty/\od b^\infty$, which corresponds to the weak competitive case.
\end{remark}

{\bf Proof of Theorem \ref{th3.1}}. The approach of this proof can be regarded as the upper and lower
solutions method. This proof will be divided into four steps. In the first one, we shall construct four positive $T$-periodic functions $\ud{U}$, $\ud{V}$, $\overline{U}$
and $\overline{V}$, for which $(\overline{U},\overline{V})$ and
$(\ud{U},\ud{V})$ can be regarded as {\it coupled ordered upper and lower
solutions} of (\ref{3.1}). In step 2, by use of the functions $\ud{U}$, $\ud{V}$, $\overline{U}$ and $\overline{V}$, we prove the existences of
$U^*,\,U_*,\,V^*$ and $V_*$. Proofs of (\ref{3.12}) and
(\ref{3.13}) will be given in the third and fourth steps, respectively.

\vskip 3pt {\bf Step 1}. The construction of $\ud{U}$, $\ud{V}$,
$\overline{U}$ and $\overline{V}$. Since $a\in{\cal C}_r(T)$, take advantage of Proposition \ref{p3.3} we know that the problem
 \bes\left\{\begin{array}{ll}
 U_t-d_1U_{xx}=U\big(a(t,x)-U\big), \ \ &0\le t\le T, \
 0<x<\yy,\\[1mm]
  B_1[U](t,0)=0,  &0\le t\le T,\\[1mm]
 U(0,x)=U(T,x), &0\le x<\yy
 \end{array}\right.\lbl{3.16} \ees
admits a unique positive solution $\overline{U}(t,x)\in{\cal C}_r(T)$, and
\bes
 &\ud a_\infty\le\dd\liminf_{x\to\infty}\frac{\overline{U}(t,x)}{x^r},\ \ \
\limsup_{x\to\infty}\frac{\overline{U}(t,x)}{x^r}\le\od a^\infty,&\lbl{3.17}\\[1mm]
&\dd\limsup_{x\to\infty}\overline{U}(t,x)\le w_2(t)\ \ \ \mbox{if} \ \
r=0&\lbl{3.18}\ees
uniformly in $[0,T]$. Moreover, $\overline{U}\le\|a\|_\infty$ by the maximum
principle. It follows that $b-h\overline{U}\in{\cal C}_r(T)$ since $b\in{\cal C}_r(T)$ and $\ud b_\infty-h\od
a^\infty>0$.
Applying Proposition \ref{p3.3} again, the problem
  \bes\left\{\begin{array}{ll}
  V_t-d_2V_{xx}=V\big(b(t,x)-h\overline{U}(t,x)-V\big), \ \
  &0\le t\le T, \ 0<x<\yy,\\[1mm]
 B_2[V](t,0)=0,  &0\le t\le T,\\[1mm]
 V(0,x)=V(T,x), &0\le x<\yy
 \end{array}\right. \lbl{3.19}\ees
has a unique positive solution $\ud{V}(t,x)\in{\cal C}_r(T)$, and
 \bes
 &\dd\ud b_\infty-h\od
 a^\infty\le\liminf_{x\to\infty}\frac{\ud{V}(t,x)}{x^r}, \ \ \
 \limsup_{x\to\infty}\frac{\ud{V}(t,x)}{x^r}\le\od b^\infty-h\ud
 a_\infty,&\lbl{3.20}\\[1mm]
 &\dd \liminf_{x\to\infty}\ud{V}(t,x)\ge z_1(t) \ \ \mbox{if} \ \
 r=0&\lbl{3.21}
 \ees
uniformly in $[0,T]$.

Similarly to the above, the problem
  \bess\left\{\begin{array}{ll}
 V_t-d_2V_{xx}=V\big(b(t,x)-V\big), \ \ &0\le t\le T, \
 0<x<\yy,\\[1mm]
 B_2[V](t,0)=0,  &0\le t\le T,\\[1mm]
 V(0,x)=V(T,x), &0\le x<\yy
 \end{array}\right.\eess
admits a unique positive solution $\overline{V}\in{\cal C}_r(T)$, and the problem
  \bess\left\{\begin{array}{ll}
  U_t-d_1U_{xx}=U\big(a(t,x)-k\overline{V}(t,x)-U\big), \ \
  &0\le t\le T, \ 0<x<\yy,\\[1mm]
 B_1[U](t,0)=0,  &0\le t\le T,\\[1mm]
 U(0,x)=U(T,x), &0\le x<\yy
 \end{array}\right. \eess
has a unique positive solution $\ud{U}\in{\cal C}_r(T)$. Moreover,
  \bes
 w_1(t)\le\liminf_{x\to\infty}\ud{U}(t,x), \ \ \limsup_{x\to\infty}\overline{V}(t,x)\le z_2(t) \ \ \mbox{if} \ \ r=0\lbl{3.23}
  \ees
uniformly in $[0,T]$.

 \vskip 2pt In addition, by Proposition \ref{p3.4} we have $\ud{U}(t,x)\le
\overline{U}(t,x)$ and $\ud{V}(t,x)\le\overline{V}(t,x)$ in $[0,T]\times[0,\yy)$, and hence in $[0,\yy)\times[0,\yy)$ as they are $T$-periodic functions in time $t$.

\vskip 2pt{\bf Step 2}. In this step we use the functions $\ud{U}$, $\ud{V}$,
$\overline{U}$ and $\overline{V}$ obtained in the above step to construct
$U^*,\,U_*,\,V^*$ and $V_*$, and prove that both
$(U^*,V_*)$ and $(U_*,V^*)$ are positive solutions of
$(\ref{3.1})$. Such a process is probably well known. For completeness, we shall provide the details.

Let $\ell>0$ and $(w_\ell,\,z_\ell)$ be the unique positive solution of the
initial-boundary value problem:
  \bes\left\{\begin{array}{ll}
 w_t-d_1 w_{xx}=w\big(a(t,x)-w-k z\big), \ \ &t>0, \ 0<x<\ell,\\[1mm]
 z_t-d_2 z_{xx}=z\big(b(t,x)-z-h w\big), \ \ &t>0, \ 0<x<\ell,\\[1mm]
 B_1[w](t,0)=B_2[z](t,0)=0,  &t\ge 0,\\[1mm]
 w(t,\ell)=\overline{U}(t,\ell), \ \ z(t,\ell)=\ud{V}(t,\ell),  &t\ge
 0,\\[1mm]
 w(0,x)=\overline{U}(0,x),\ \ z(0,x)=\ud{V}(0,x), \ \
  &0\le x\le \ell.
 \end{array}\right. \lbl{3.24}\ees
Since $\overline{U}$ and $\ud{V}$ are positive functions and satisfy
 \bess\left\{\begin{array}{ll}
 \overline{U}_t-d_1\overline{U}_{xx}>\overline{U}\big(a(t,x)-\overline{U}-k\ud{V}\big),
 \ \ &t>0, \ 0<x<\ell,\\[1mm]
 \ud{V}_t-d_2\ud{V}_{xx}=\ud{V}\big(b(t,x)-\ud{V}-h\overline{U}\big),
 \ \ &t>0, \ 0<x<\ell,\\[1mm]
 B_1[\overline{U}](t,0)=B_2[\ud{V}](t,0)=0,  &t\ge 0,
 \end{array}\right. \eess
one can use Lemma \ref{l3.1} to deduce $w_\ell\le \overline{U},\,z_\ell\ge
\ud{V}$ in $[0,\infty)\times[0,\ell]$. For the non-negative integer $n$,
we define
 \[w_\ell^n(t,x)=w_\ell(t+nT,x), \ \ \ z_\ell^n(t,x)=z_\ell(t+nT,x).\]
Note that $a(t,x),\,b(t,x),\,\overline{U}(t,x)$ and $\ud{V}(t,x)$ are
$T$-periodic functions in time $t$, it follows that $(w_\ell^n,z_\ell^n)$ satisfies
 \bess\left\{\begin{array}{ll}
 (w^n_\ell)_t-d_1(w^n_\ell)_{xx}=w^n\big(a(t,x)-w^n-k z^n\big), \ &0<t\le T, \
 0<x<\ell,\\[1mm]
 (z^n_\ell)_t-d_2 (z^n_\ell)_{xx}=z^n\big(b(t,x)-z^n-h w^n\big), \ &0<t\le T, \
 0<x<\ell,\\[1mm]
 B_1[w^n](t,0)=B_2[z^n](t,0)=0,  &0\ge t\le T,\\[1mm]
 w^n_\ell(t,\ell)=\overline{U}(t,\ell), \ \ z^n_\ell(t,\ell)=\ud{V}(t,\ell),
 &0\ge t\le T,\\[1mm]
 w^n_\ell(0,x)=w_\ell(nT,x),\ \ z^n_\ell(0,x)=z_\ell(nT,x), \ \  &0\le x\le \ell.
 \end{array}\right.\eess
Remember
 \bess
 &w^1_\ell(0,x)=w_\ell(T,x)\le\overline{U}(T,x)=\overline{U}(0,x)=w_\ell(0,x),&
 \\[1mm]
 &z^1_\ell(0,x)=z_\ell(T,x)\ge \ud{V}(T,x)=\ud{V}(0,x)=z_\ell(0,x),&\eess
it is derived by Lemma \ref{l3.1} that $w^1_\ell\le w_\ell$, $z^1_\ell\ge z_\ell$ in
$[0,T]\times[0,\ell]$. And then
 \[w^2_\ell(0,x)=w^1_\ell(T,x)\le w_\ell(T,x)=w^1_\ell(0,x), \ \
 z^2_\ell(0,x)=z^1_\ell(T,x)\ge
  z_\ell(T,x)=z^1_\ell(0,x).\]
Apply Lemma \ref{l3.1} once again, we have $w_\ell^2\le w_\ell^1$, $z_\ell^2\ge
z_\ell^1$ in $[0,T]\times[0,\ell]$. Utilizing the inductive method we can show that
$w_\ell^n$ and  $z_\ell^n$ are, respectively, decreasing and increasing in $n$. So,
there exist two non-negative functions $\overline{U}_\ell$, $\ud{V}_\ell$ such
that $w_\ell^n\to\overline{U}_\ell$, $z_\ell^n\to\ud{V}_\ell$ pointwise in
$[0,T]\times[0,\ell]$ as $n\to\infty$.  Obviously,
$\overline{U}_\ell(0,x)=\overline{U}_\ell(T,x)$,
$\ud{V}_\ell(0,x)=\ud{V}_\ell(T,x)$ as $w^{n+1}_\ell(0,x)=w^n_\ell(T,x)$,
$z^{n+1}_\ell(0,x)=z^n_\ell(T,x)$. Based on the regularity theory for
parabolic equations and compact argument, it can be proved that there exists a
subsequence $\{n_i\}$, such that $w_\ell^{n_i}\to\overline{U}_\ell$,
$z_\ell^{n_i}\to\ud{V}_\ell$ in $C^{1,2}([0,T]\times[0,\ell])$ as $i\to\infty$,
and $(\overline{U}_\ell,\,\ud{V}_\ell)$ satisfies the first four equations of
(\ref{3.24}). Therefore, $(\overline{U}_\ell,\,\ud{V}_\ell)$ satisfies
 \bes\left\{\begin{array}{ll}
 U_t-d_1 U_{xx}=U\big(a(t,x)-U-k V\big), \ \ &0\le t\le T,
 \ 0<x<\ell,\\[1mm]
 V_t-d_2 V_{xx}=V\big(b(t,x)-V-h U\big), \ \ &0\le t\le T,
 \ 0<x<\ell,\\[1mm]
 B_1[U](t,0)=B_2[V](t,0)=0,  &0\le t\le T,\\[1mm]
 U(t,\ell)=\overline{U}(t,\ell), \ \ V(t,\ell)=\ud{V}(t,\ell),  &0\le
 t\le T,\\[1mm]
 U(0,x)=U(T,x),\ \ V(0,x)=V(T,x), \ \
  &0\le x\le \ell.
 \end{array}\right. \lbl{3.25}\ees
Evidently, $\overline{U}_\ell(t,x)>0,\,\ud{V}_\ell(t,x)>0$ in $(0,T]\times[0,\ell]$
because $\overline{U}(t,\ell)>0,\,\ud{V}(t,\ell)>0$ in $(0,T]$. This shows that
$(\overline{U}_\ell,\,\ud{V}_\ell)$ is a positive solution of (\ref{3.25}).

Let $(\varphi_\ell,\,\psi_\ell)$ be the unique positive solution of the following
initial-boundary value problem
  \bess\left\{\begin{array}{ll}
 \varphi_t-d_1 \varphi_{xx}=\varphi\big(a(t,x)-\varphi-k \psi\big), \ &t>0, \
 0<x<\ell,\\[1mm]
 \psi_t-d_2 \psi_{xx}=\psi\big(b(t,x)-\psi-h \varphi\big), \ &t>0, \
 0<x<\ell,\\[1mm]
 B_1[\varphi](t,0)=B_2[\psi](t,0)=0,  &t\ge 0,\\[1mm]
 \varphi(t,\ell)=\ud{U}(t,\ell), \ \ \psi(t,\ell)=\overline{V}(t,\ell),  &t\ge
 0,\\[1mm]
 \varphi(0,x)=\ud{U}(0,x),\ \ \psi(0,x)=\overline{V}(0,x), \ \
  &0\le x\le \ell.
 \end{array}\right.\eess
Since $\ud{U}\le\overline{U}$, $\ud{V}\le\overline{V}$ in $[0,\yy)\times[0,\yy)$, we have $w_\ell\ge\varphi_\ell\ge\ud{U},\,z_\ell\le\psi_\ell\le\overline{V}$ in
$[0,\yy)\times[0,\ell]$ by Lemma \ref{l3.1}. Similarly to the above, there exist two
positive functions $\ud{U}_\ell$, $\overline{V}_\ell$, and a subsequence
$\{n_i\}$, such that $\varphi_\ell(t+n_iT,x)\to\ud{U}_\ell(t,x)$,
$\psi_\ell(t+n_iT,x)\to\overline{V}_\ell(t,x)$ in $C^{1,2}([0,T]\times[0,\ell])$ as
$i\to\infty$, and $(\ud{U}_\ell,\,\overline{V}_\ell)$ solves
 \bes\left\{\begin{array}{ll}
 U_t-d_1 U_{xx}=U\big(a(t,x)-U-k V\big), \ \ &0\le t\le T,
 \ 0<x<\ell,\\[1mm]
 V_t-d_2 V_{xx}=V\big(b(t,x)-V-h U\big), \ \ &0\le t\le T,
 \ 0<x<\ell,\\[1mm]
 B_1[U](t,0)=B_2[V](t,0)=0,  &0\le t\le T,\\[1mm]
 U(t,\ell)=\ud{U}(t,\ell), \ \ V(t,\ell)=\overline{V}(t,\ell),  &0\le
 t\le T,\\[1mm]
 U(0,x)=U(T,x),\ \ V(0,x)=V(T,x), \ \
  &0\le x\le \ell.
 \end{array}\right. \lbl{3.26}\ees

Recall that $\ud{U}\le\varphi_\ell\le w_\ell\le\overline{U},\,\ud{V}\le
z_\ell\le\psi_\ell\le\overline{V}$ in $[0,\yy)\times[0,\ell]$, it is immediately
to get
 \[\ud{U}\le\ud{U}_\ell\le\overline{U}_\ell\le\overline{U}, \ \
 \,\ud{V}\le\ud{V}_\ell\le\overline{V}_\ell\le\overline{V} \ \ \ \mbox{in}
 \ \ [0,T]\times[0,\ell]\]
for any given $\ell>0$. By use of the regularity theory for parabolic equations and
compact argument, we can show that there exist two subsequences
$\{(\overline{U}_{\ell_j},\ud{V}_{\ell_j})\}$,
$\{(\ud{U}_{\ell_j},\overline{V}_{\ell_j})\}$ and four positive functions
$U^*,\,U_*,\,V^*$, $V_*$, such that, for any $L>0$,
 \bess
 (\overline{U}_{\ell_j},\ud{V}_{\ell_j})\to (U^*,V_*), \ \
 (\ud{U}_{\ell_j},\overline{V}_{\ell_j})\to(U_*,V^*)
 \ \ \mbox{in}\ \ \big[C^{1,2}([0,T]\times[0,L])\big]^2\eess
as $j\to\infty$. Obviously,
 \bes
 \ud{U}\le U_*\le U^*\le\overline{U}, \ \ \,\ud{V}\le
 V_*\le V^*\le\overline{V} \ \ \ \mbox{in} \ \
 [0,T]\times[0,\infty).\lbl{3.27}\ees
Remember (\ref{3.25}) and (\ref{3.26}), it is easy to see that both
$(U^*,V_*)$ and $(U_*,V^*)$ are positive $T$-periodic solutions
of (\ref{3.1}).

\vskip 2pt{\bf Step 3}. Prove (\ref{3.12}). Let $(U,V)$ be a positive
solution of (\ref{3.1}). We only prove $U\le U^*,\,V\ge V_*$ as
the proof of $U\ge U_*,\,V\le V^*$ is similarly.

Because $U$ satisfies
 \bess\left\{\begin{array}{ll}
 U_t-d_1 U_{xx}<U\big(a(t,x)-U\big), \ \ &0\le t\le T, \
 0<x<\yy,\\[1mm]
  B_1[U](t,0)=0,  &0\le t\le T,\\[1mm]
 U(0,x)=U(T,x), &0\le x<\yy,
 \end{array}\right. \eess
we have $U(t,x)\le\|a\|_\infty$ for all $(t,x)\in[0,T]\times[0,\infty)$.
Choose $\ell\gg 1$. Using Lemma \ref{l3.2} and \cite[Theorem 28.1]{Hess},
respectively, we have that the following problems
 \bess\left\{\begin{array}{ll}
 \zeta_t-d_1\zeta_{xx}=\zeta\big(a(t,x)-\zeta\big), \ \ &0\le
 t\le T, \ 0<x<\ell,\\[1mm]
 B_1[\zeta](t,0)=0, \ \ \zeta(t,\ell)=\|a\|_\infty, \ \  &0\le t\le
 T,\\[1mm]
 \zeta(0,x)=\zeta(T,x), &0\le x\le \ell
 \end{array}\right.\eess
and
 \bess\left\{\begin{array}{ll}
 \zeta_t-d_1\zeta_{xx}=\zeta\big(a(t,x)-\zeta\big), \ \ &0\le
 t\le T, \ 0<x<\ell,\\[1mm]
 B_1[\zeta](t,0)=0, \ \ \zeta(t,\ell)=0, \ \  &0\le t\le T,\\[1mm]
 \zeta(0,x)=\zeta(T,x), &0\le x\le \ell
 \end{array}\right.\eess
have unique positive solutions $\overline{\zeta}_\ell$ and $\ud{\zeta}_\ell$.
We can apply Lemma \ref{l3.3} to conclude that  $U\le\overline{\zeta}_\ell$ and $\ud{\zeta}_\ell\le\overline{\zeta}_\ell$ in $[0,T]\times[0,\ell]$,
$\overline{\zeta}_\ell$ and $\ud{\zeta}_\ell$ are decreasing and increasing in
$\ell$, respectively. Similarly to the proof of Proposition \ref{p3.4}, we have
$\dd\lim_{\ell\to\infty}\overline{\zeta}_\ell=\overline{U}$ in $C^{1, 2}([0,T]\times[0,L])$ for any $L>0$ since $\overline{U}$ is the unique positive
solution of (\ref{3.16}). Therefore, $U\le\overline{U}$ in
$[0,T]\times[0,\infty)$. Similarly, we have $V\ge\ud{V}$ in
$[0,T]\times[0,\infty)$. It follows that $(U,V)$ satisfies
  \bes\left\{\begin{array}{ll}
 U_t-d_1 U_{xx}=U\big(a(t,x)-U-k V\big), \ \ &0<t\le T, \
 0<x<\ell,\\[1mm]
 V_t-d_2 V_{xx}=V\big(b(t,x)-V-h U\big), \ \ &0<t\le T, \
 0<x<\ell,\\[1mm]
 B_1[U](t,0)=B_2[V](t,0)=0,  &0\le t\le T,\\[1mm]
 U(t,\ell)\le\overline{U}(t,\ell), \ \ V(t,\ell)\ge\ud{V}(t,\ell),
 &0\le t\le T,\\[1mm]
 U(0,x)\le\overline{U}(0,x),\ \ V(0,x)\ge\ud{V}(0,x), \ \
  &0\le x\le \ell.
 \end{array}\right. \lbl{3.28}\ees
Applying Lemma \ref{l3.1} to the problems (\ref{3.24}) and (\ref{3.28}), we get that
  \bes
  U\le w_\ell, \ \ V\ge z_\ell \ \ \ \mbox{in}\ \
  [0,T]\times[0,\ell].\lbl{3.29}\ees
It can be seen from the arguments of step 2 that
 \bes\lim_{j\to\infty}\dd\lim_{n\to\infty}w_{\ell_j}(t+nt,x)=U^*, \ \ \
 \lim_{j\to\infty}\lim_{n\to\infty}z_{\ell_j}(t+nt,x)=V_*\lbl{3.30}\ees
in $C^{1,2}([0,T]\times[0,L])$ for any $L>0$. According to
$U(t+nT,x)=U(t,x)$ and $V(t+nT,x)=V(t,x)$, it is derived from
(\ref{3.29}) and (\ref{3.30}) that $U\le U^*,\,V\ge V_*$.

\vskip 2pt{\bf Step 4}. When $r=0$ in ${\bf(H1)}$. Summarizing (\ref{3.18}),
(\ref{3.21}) and (\ref{3.23}) we obtain
 \bess\begin{array}{ll}
 w_1(t)\le\dd\liminf_{x\to\infty}\ud{U}(t,x),\ \
\limsup_{x\to\infty}\overline{U}(t,x)\le w_2(t),\\
z_1(t)\le\dd\liminf_{x\to\infty}\ud{V}(t,x),\ \
\limsup_{x\to\infty}\overline{V}(t,x)\le z_2(t)
 \end{array}\eess
uniformly in $[0,T]$. Combining these facts with (\ref{3.27}) and (\ref{3.12}), we can derive (\ref{3.13}). The proof is complete. \ \ \ \fbox{}

\section{Spreading-vanishing dichotomy of the problem (\ref{1.2})}
\setcounter{equation}{0}

We first state a lemma, by which the vanishing phenomenon is immediately obtained.
Moreover, this lemma will play an important role in the establishment of criteria
for spreading and vanishing.
  \begin{lem}\hspace{-1mm}{\rm(\cite[Lemma 3.1]{Wjde15})}\lbl{l4.1}\, Let $d,\mu$ and $B$ be as above, $C\in\mR$. Assume
that functions $g\in C^1([0,\infty))$, $\vp\in
C^{\frac{1+\nu}2,1+\nu}([0,\infty)\times[0,g(t)])$ and satisfy $g(t)>0$,
$\vp(t,x)>0$ for $t\ge 0$ and $0<x<g(t)$. We further suppose that
$\dd\lim_{t\to\infty} g(t)<\infty$, $\dd\lim_{t\to\infty} g'(t)=0$ and there exists a
constant $K>0$ such that $\|\vp(t,\cdot)\|_{C^1[0,\,g(t)]}\leq K$ for $t>1$. If
$(\vp,g)$ satisfies
  \bess\left\{\begin{array}{lll}
 \vp_t-d \vp_{xx}\geq C\vp, &t>0, \ 0<x<g(t),\\[1mm]
 B[\vp]=0, \ &t\ge 0, \ x=0,\\[1mm]
 \vp=0,\ \ g'(t)\geq-\mu \vp_x, \ &t\ge 0, \ x=g(t),
 \end{array}\right.\eess
then $\dd\lim_{t\to\infty}\max_{0\leq x\leq g(t)}\vp(t,x)=0$.
 \end{lem}

Applying (\ref{2.3}) and Lemma \ref{l4.1}, we have the following
result.
 \begin{theo}\lbl{th4.1}$($Vanishing$)$\, If $s_\infty<\infty$, then
 \bes
 \lim_{t\to\infty}\|u(t,\cdot)\|_{C([0,s(t)])}=0, \ \ \
 \lim_{t\to\infty}\|v(t,\cdot)\|_{C([0,s(t)])}=0.\lbl{4.1}
 \ees
This shows that if the two species cannot spread successfully, they will extinct in
the long run.
\end{theo}

For any given $\ell>0$, let $\lm_1(\ell;d_1,a)$ and $\gamma_1(\ell;d_2,b)$ be the
principal eigenvalues of the $T$-periodic eigenvalue problems
 \bes\left\{\begin{array}{ll}
 \phi_t-d_1\phi_{xx}-a(t,x)\phi=\lm\phi, \ \ &0\le t\le T, \ 0<x< \ell,\\[1mm]
 B_1[\phi](t,0)=0, \ \ \phi(t, \ell)=0,\ \ &0\le t\le T, \\[1mm]
 \phi(0,x)=\phi(T,x),&0\le x\le \ell
 \end{array}\right.\lbl{3.9}\ees
and
 \bes\left\{\begin{array}{ll}
 \psi_t-d_2\psi_{xx}-b(t,x)\psi=\gamma\psi, \ \ &0\le t\le T, \ 0<x< \ell,\\[1mm]
 B_2[\psi](t,0)=0, \ \ \psi(t, \ell)=0,\ \ &0\le t\le T, \\[1mm]
 \psi(0,x)=\psi(T,x),&0\le x\le \ell,
 \end{array}\right.\lbl{3.10}\ees
respectively.
  \begin{theo}\lbl{th4.2}$($Spreading$)$\, Assume that {\bf(H1)} holds and
$s_\infty=\infty$. If the condition $(\ref{3.11})$ is true, then we have
 \bes
 & U_*(t,x)\le\dd\liminf_{n\to\infty}u(t+nT,x), \ \
 \limsup_{n\to\infty}u(t+nT,x)\leq U^*(t,x),&\lbl{4.2}\\
 &V_*(t,x)\le\dd\liminf_{n\to\infty}v(t+nT,x), \ \
 \limsup_{n\to\infty}v(t+nT,x)\leq V^*(t,x)&\lbl{4.3}
 \ees
uniformly in $[0, T]\times[0,L]$ for any $L>0$, where $(U^*,V_*)$ and
$(U_*,V^*)$ are positive $T$-periodic solutions of $(\ref{3.1})$ obtained
in Theorem $\ref{th3.1}$.
\end{theo}

{\bf Proof}. The method used here is an {\it iterative process}. This proof not only
gives the long time behavior of $(u,v)$ but also the existence of
$(U^*,V_*)$ and $(U_*,V^*)$.

The proof is divided into four steps. Let functions
$\overline{U},\,\ud{U},\,\overline{V}$ and $\ud{V}$ be given in the proof
of Theorem \ref{l3.1}. In the first two steps we shall prove, respectively, that
  \bes
  \limsup_{n\to\infty}u(t+nT,x)\le\overline{U}(t,x)\ \ \mbox{uniformly in}
  \ [0,T]\times[0,L]\lbl{4.4}\ees
and
 \bes
\liminf_{n\to\infty}v(t+nT,x)\ge\ud{V}(t,x)\ \ \mbox{uniformly in} \
[0,T]\times[0,L]
 \lbl{4.4a}\ees
for any given $L>0$. In the third one, we shall construct four sequences $\{\overline{U}_i\}$,
$\{\ud{V}_i\}$, $\{\ud{U}_i\}$ and $\{\overline{V}_i\}$ satisfying
 \bes\ud{U}\le\ud{U}_1\le\cdots\le\ud{U}_i\le\overline{U}_i\le\cdots\overline{U}_1\le\overline{U},\
 \ \
 \ud{V}\le\ud{V}_1\le\cdots\le\ud{V}_i\le\overline{V}_i\le\cdots\overline{V}_1\le\overline{V}.
  \lbl{4.4b}\ees
Proofs of (\ref{4.2}) and (\ref{4.3}) will be given in the last step.

{\bf Step 1}. \,Define
 \[\phi(x)=\left\{\begin{array}{ll}
  u_0(x), \ \ & 0\leq x\leq s_0,\\[1mm]
  0, \ \ & x\geq s_0,
  \end{array}\right.\]
and let $w(t,x)$ be the unique positive solution of
 \bess\left\{\begin{array}{lll}
 w_t-d_1w_{xx}=w\big(a(t,x)-w\big),\ &t>0,\ 0<x<\yy,\\[1mm]
 B_1[w](t,0)=0,\ \ \ &t>0,\\[1mm]
  w(0,x)=\phi(x), &0\le x<\yy.
  \end{array}\right.\eess
In view of Proposition \ref{p3.3}, it follows that $\dd\lim_{n\to\infty}w(t+nT,x)=\overline{U}(t,x)$
uniformly in $[0,T]\times[0, L]$, where $\overline{U}(t,x)$ is
the unique positive solution of (\ref{3.16}). On the other hand, by the comparison principle, we have $u(t,x)\leq w(t,x)$ for all $t>0$ and $0\leq x\leq
s(t)$. Thanks to $s_\infty=\infty$, we get
(\ref{4.4}).

{\bf Step 2}.\, For any $\varepsilon>0$, denote
$b_\varepsilon(t,x)=b(t,x)-h\big(\overline{U}(t,x)+\varepsilon(1+x)^r\big)$. It
follows from (\ref{1.4}) and (\ref{3.17}) that
 \bess
  \ud b_\infty-h(\od a^\infty+\varepsilon)
  \le\liminf_{x\to\infty}\frac{b_\varepsilon(t,x)}{x^r}
  \le\limsup_{x\to\infty}\frac{b_\varepsilon(t,x)}{x^r}\le\od b^\infty-h(\ud
  a_\infty+\varepsilon)\eess
uniformly in $[0,T]$. Since $\ud b_\infty>h\od a^\infty$, there exists
$\varepsilon_0>0$ such that $\ud b_\infty>h(\od a^\infty+\varepsilon)$, and hence
$b_\ep\in {\cal C}_r(T)$ for all $0<\varepsilon\le\varepsilon_0$. For such fixed
$\varepsilon$, by Proposition \ref{p3.2}, there exists $\ell_0^\varepsilon>L$ such
that $\gamma_1(\ell;d_2,b_\varepsilon)<0$ for all $\ell\ge \ell_0^\varepsilon$.

For any fixed $0<\varepsilon<\varepsilon_0$ and $\ell>\ell_0^\varepsilon$, capitalize on (\ref{4.4}) and $s_\yy=\yy$, there exists $\tau\gg 1$ such that
 \[s(t)>\ell, \ \ u(t,x)<\overline{U}(t,x)+\varepsilon(1+x)^r\ff t\ge\tau, \ 0\leq
 x\leq \ell.\]
Consider the following auxiliary $T$-periodic boundary value problem
\bess\left\{\begin{array}{ll}
 Z_t-d_2 Z_{xx}=Z\big(b_\varepsilon(t,x)-Z\big),\ \ &0\le t\le
 T, \ 0<x<\ell,\\[1mm]
  B_2[Z](t,0)=Z(t,\ell)=0, &0\le t\le T,\\[1mm]
 Z(0,x)=Z(T,x), &0\le x\le \ell.\end{array}\right.
 \eess
Since $\gamma_1(\ell;d_2,b_\varepsilon)<0$, utilizing Theorem 28.1 of \cite{Hess},
the above problem admits a unique positive solution, denoted by $Z^\varepsilon_\ell(t,x)$. Let
$V^\varepsilon_\ell(t,x)$ be the unique positive solution of the following initial-boundary value problem
 \bess\left\{\begin{array}{lll}
 V_t-d_2 V_{xx}=V\big(b_\varepsilon(t,x)-V\big),\ \ &t>\tau, \
 0<x<\ell,\\[1mm]
 B_2[V](t,0)=0, \ \ V(t,\ell)=0,\ \ \ &t\ge\tau,\\[1mm]
 V(\tau,x)=\sigma Z^\ep_\ell(\tau,x),&x\in [0,\ell],
 \end{array}\right.\eess
where $0<\sigma<1$ is so small that $\sigma
Z^\varepsilon_\ell(\tau,x)<v(\tau,x)$ in $[0,\ell]$. Obviously, the function
$\chi:=\sigma Z^\varepsilon_\ell$ satisfies
 \bess\left\{\begin{array}{lll}
 \chi_t-d_2 \chi_{xx}<\chi\big(b_\varepsilon(t,x)-\chi\big),\ \ &t>\tau, \
 0<x<\ell,\\[1mm]
 B_2[\chi](t,0)=0, \ \ \chi(t,\ell)=0,\ \ \ &t\ge\tau,\\[1mm]
 \chi(\tau,x)=\sigma Z^\ep_\ell(\tau,x),&x\in [0,\ell].
 \end{array}\right.\eess
By the comparison principle,
\bess
 v(t,x)\geq V^\varepsilon_\ell(t,x)\ge\chi(t,x)\ff t\geq\tau, \ 0\leq x\leq \ell.\eess

Using the arguments of step 2 in the proof of Theorem \ref{th3.1}, we can prove that
$\dd\lim_{n\to\infty}V^\varepsilon_\ell(t+nT,x)=Z^\varepsilon_\ell(t,x)$ in
$C^{1,2}([0,T]\times[0,\ell])$ and
$\dd\lim_{\ell\to\infty}Z^\varepsilon_\ell(t,x)=Z^\varepsilon(t,x)$ in $C^{1,
2}([0,T]\times[0,L])$, where $Z^\varepsilon$ is the
unique positive solution of $T$-periodic boundary value problem
 \bess\left\{\begin{array}{ll}
 Z_t-d_2 Z_{xx}=Z\big(b_\varepsilon(t,x)-Z\big),\ \ &0\le t\le
 T, \ 0<x<\yy,\\[1mm]
  B_2[Z](t,0)=Z(t,\ell)=0, &0\le t\le T,\\[1mm]
 Z(0,x)=Z(T,x), &0\le x<\yy.\end{array}\right.
 \eess
The existence and uniqueness of $Z^\varepsilon$ is guaranteed by Proposition
\ref{p3.3}. It follows that
 \[\liminf_{n\to\infty}v(t+nT,x)\ge Z^\varepsilon(t,x)\ \ \mbox{
 uniformly for} \ (t,x)\in[0,T]\times[0,L].\]
Note that $b_\varepsilon(t,x)\to b(t,x)-h\overline{U}(t,x)$ as $\ep\to 0$ and $\ud{V}(t,x)$ is the unique positive solution of (\ref{3.19}), by the continuous dependence of solution with respect to parameter, we have that
$\dd\lim_{\varepsilon\to 0} Z^\varepsilon(t,x)=\ud{V}(t,x)$ uniformly in $[0,T]\times[0,L]$. Thus, (\ref{4.4a}) holds.

{\bf Step 3}.\, In view of (\ref{1.4}), (\ref{3.11}) and (\ref{3.20}), we see that $a-k\ud{V}\in{\cal C}_r(T)$. Same as the second step, it can be deduced that
$\dd\limsup_{n\to\infty}u(t+nT,x)\le\overline{U}_1(t,x)$ locally uniformly for
$(t,x)\in[0,T]\times[0,\infty)$, where $\overline{U}_1$ is the unique positive solution of
 \bess\left\{\begin{array}{ll}
 U_t-d_1U_{xx}=U\big(a(t,x)-k\ud{V}(t,x)-U\big), \ \ &0\le
 t\le T, \ 0<x<\yy,\\[1mm]
  B_1[U](t,0)=0,  &0\le t\le T,\\[1mm]
 U(0,x)=U(T,x), &0\le x<\yy.
 \end{array}\right. \eess
Similarly, $\dd\liminf_{n\to\infty}v(t+nT,x)\ge\ud{V}_1(t,x)$ locally uniformly in $[0,T]\times[0,\infty)$, where $\ud{V}_1$ is the unique positive solution of
 \bess\left\{\begin{array}{ll}
  V_t-d_2V_{xx}=V\big(b(t,x)-h\overline{U}_1(t,x)-V\big), \ \
  &0\le t\le T, \ 0<x<\yy,\\[1mm]
 B_2[V](t,0)=0,  &0\le t\le T,\\[1mm]
 V(0,x)=V(T,x), &0\le x<\yy.
 \end{array}\right. \eess

Repeating the above procedure, we can find two sequences $\{\overline{U}_i\}$ and
$\{\ud{V}_i\}$ such that
  \bes\limsup_{n\to\infty}u(t+nT,x)\le\overline{U}_i(t,x), \ \
  \liminf_{n\to\infty}v(t+nT,x)\ge\ud{V}_i(t,x)\lbl{4.5}\ees
locally uniformly for $(t,x)\in[0,T]\times[0,\infty)$, here $\overline{U}_i$ and
$\ud{V}_i$ are the unique positive solutions of
  \bes\left\{\begin{array}{ll}
 U_t-d_1U_{xx}=U\big(a(t,x)-k\ud{V}_{i-1}(t,x)-U\big), \ \
 &0\le t\le T, \ 0<x<\yy,\\[1mm]
  B_1[U](t,0)=0,  &0\le t\le T,\\[1mm]
 U(0,x)=U(T,x), &0\le x<\yy
 \end{array}\right.\lbl{4.6} \ees
and
 \bes\left\{\begin{array}{ll}
  V_t-d_2V_{xx}=V\big(b(t,x)-h\overline{U}_i(t,x)-V\big), \ \
  &0\le t\le T, \ 0<x<\yy,\\[1mm]
 B_2[V](t,0)=0,  &0\le t\le T,\\[1mm]
 V(0,x)=V(T,x), &0\le x<\yy,
 \end{array}\right. \lbl{4.7}\ees
respectively.

In the same way we can get two sequences $\{\ud{U}_i\}$ and
$\{\overline{V}_i\}$ such that
  \bes\liminf_{n\to\infty}u(t+nT,x)\ge\ud{U}_i(t,x), \ \
  \limsup_{n\to\infty}v(t+nT,x)\le\overline{V}_i(t,x)\lbl{4.8}\ees
locally uniformly in $[0,T]\times[0,\infty)$, here $\overline{V}_i$ and
$\ud{U}_i$ are the unique positive solutions of
  \bes\left\{\begin{array}{ll}
  V_t-d_2V_{xx}=V\big(b(t,x)-h\ud{U}_{i-1}(t,x)-V\big), \ \
  &0\le t\le T, \ 0<x<\yy,\\[1mm]
 B_2[V](t,0)=0,  &0\le t\le T,\\[1mm]
 V(0,x)=V(T,x), &0\le x<\yy
 \end{array}\right.\lbl{4.9}\ees
and
  \bes\left\{\begin{array}{ll}
 U_t-d_1U_{xx}=U\big(a(t,x)-k\overline{V}_i(t,x)-U\big), \ \
 &0\le t\le T, \ 0<x<\yy,\\[1mm]
  B_1[U](t,0)=0,  &0\le t\le T,\\[1mm]
 U(0,x)=U(T,x), &0\le x<\yy,
 \end{array}\right. \lbl{4.10}\ees
respectively.

Apply Proposition \ref{p3.4}, we can show that (\ref{4.4b}) holds.

{\bf Step 4}.\, Now we prove (\ref{4.2}) and (\ref{4.3}). Remembering (\ref{4.4b}),
make use of the regularity theory for parabolic equations and compact argument, we assert that there exist four positive $T$-periodic functions
$U^\yy,\,U_\yy,\,V^\yy,\,V_\yy\in C^{1+\frac{\nu}2,
2+\nu}([0,T]\times[0,\infty))$, such that
 \[(\overline{U}_i,\ud{U}_i,\overline{V}_i,\ud{V}_i)\to
 (U^{\yy},\,U_\yy,\,V^\yy,\,V_\yy) \ \ \ \mbox{as} \ \
 i\to\infty\]
in $C^{1, 2}([0,T]\times[0,K])$ for any $K>0$. Taking $i\to\infty$ in (\ref{4.6}),
(\ref{4.7}), (\ref{4.9}) and (\ref{4.10}), it derives that both
$(U^{\yy},\,V_\yy)$ and $(U_\yy,\,V^\yy)$ are positive solutions
of (\ref{3.1}). Hence, by (\ref{3.12}),
 \bes
 U_*\le U_\yy\le U^\yy\le U^*, \ \ \ V_*\le V_\yy\le
 V^\yy\le V^*.
 \lbl{4.11}\ees

Arguing as step 3 in the proof of Theorem \ref{th3.1}, it can be shown that any
positive solution $(U,V)$ of (\ref{3.1}) must satisfy $\ud{U}_i\le
U\le\overline{U}_i$, $\ud{V}_i\le V\le\overline{V}_i$ for all $i$.
Thus $U_\yy\le U\le U^\yy$ and $V_\yy\le V\le V^\yy$.
Since $(U^*,V_*)$ and $(U_*,V^*)$ are positive solutions of
$(\ref{3.1})$, we have
 \[U_\yy\le U_*\le U^*\le U^\yy, \ \ \ V_\yy\le V_*\le
 V^*\le V^\yy.\]
Recall (\ref{4.11}), it yields that $U_\yy=U_*,\,U^*=U^\yy,\,V_\yy=V_*,\,V^*=
V^\yy$. Letting $i\to\infty$ in (\ref{4.5}) and (\ref{4.8}), the required
results (\ref{4.2}) and (\ref{4.3}) are obtained.
 \ \ \ \fbox{}

\section{Criteria for spreading and vanishing of $(\ref{1.2})$}
\setcounter{equation}{0}{\setlength\arraycolsep{2pt}

Throughout this section, we assume that $(u,v,s)$ is the unique solution of
$(\ref{1.2})$. We first state a comparison principle.
 \begin{lem}$($Comparison principle$)$\label{l5.1}\, Let $\tau>0$, $\bar s\in
C^1([0,\tau])$ and $\bar s(t)>0$ in $[0,\tau]$. Let $\bar u, \bar v\in
C(\overline{O})\bigcap C^{1,2}(O)$ with $O=\{(t,x): 0<t\le\tau,\, 0<x<\bar s(t)\}$.
Assume that $(\bar u,\bar v,
\bar s)$ satisfies
\bess
 \left\{\begin{array}{ll}
  \bar u_t-d_1\bar u_{xx}\geq\bar u\big(a(t,x)-\bar u\big),\ \ &0<t\le\tau,\
  0<x<\bar s(t),\\[1mm]
 \bar v_t-d_2\bar v_{xx}\geq \bar v\big(b(t,x)-\bar v\big),&0<t\le\tau,\ 0<x<\bar
 s(t),\\[1mm]
 B_1[\bar u](t,0)\geq 0,\ \ B_2[\bar v](t,0)\geq 0,&0\le t\le\tau,\\[1mm]
 \bar u(t,\bar s(t))=\bar v(t,\bar s(t))=0,&0\le t\le\tau,\\[1mm]
 \bar s'(t)\geq-\mu[\bar u_x(t,\bar s(t))+\rho \bar v_x(t,\bar s(t))],\ \ &0\le
 t\le\tau.
 \end{array}\right.
\eess
If $\bar s(0)\geq s_0,\ \bar u(0,x)\geq 0,\ \bar v(0,x)\geq 0$ in $[0,\bar s(0)]$,
and
 $u_0(x)\leq\bar u(0,x),\ v_0(x)\leq\bar v(0,x)$ in $[0,s_0]$,
then the solution $(u,v,s)$ of {\rm(\ref{1.2})} satisfies
\[s(t)\leq\bar s(t)\ \ {\rm in}\,\ [0,\tau];\ \ \ u(t,x)\leq\bar u(t,x), \
v(t,x)\leq\bar v(t,x)\ \ {\rm in}\,\ \overline Q,\]
where $Q=\{(t,x):\, 0<t\le\tau,\ 0<x<s(t)\}$.
\end{lem}

{\bf Proof}.\, The proof is same as that of \cite[Lemma 4.1]{Wjde14} (see also the
argument of \cite[Lemma 5.1]{GW}), we omit the details. \ \ \ $\Box$

\vskip 2pt Define
 \[{\cal E}=\{\ell>0:\, \lm_1(\ell;d_1,a)=0\,\ {\rm or}\,\ \gamma_1(\ell;d_2,b)=0\}.\]
If one of the functions $a(t,x)$ and $b(t,x)$ satisfies the condition {\bf(A)}, make use of Proposition \ref{p3.2}, it yields that ${\cal E}\not=\emptyset$. Especially, when
the assumption {\bf(H1)} holds, both $a(t,x)$ and $b(t,x)$ satisfy the condition {\bf(A)}, and hence ${\cal E}\not=\emptyset$.

Now we give a necessary condition of vanishing.
 \begin{lem}\lbl{l5.2}\, Assume that ${\cal E}\not=\emptyset$ and set $s^*=\min{\cal E}$. If $s_\infty<\infty$, then $s_\infty\leq s^*$. Hence, $s_0\ge s^*$ implies
$s_\infty=\infty$ for all $\mu>0$.
\end{lem}

{\bf Proof}.\, First of all, $s^*>0$ since $\lm_1(\ell;d_1,a)>0$,
$\gamma_1(\ell;d_2,b)>0$ when $0<\ell\ll 1$. Without loss of generality we assume
that $\lm_1(s^*;d_1,a)=0$.

If $s_\infty>s^*$, then $\lm_1(s_\infty;d_1,a)<0$ since $\lm_1(\ell;d_1,a)$ is
strictly decreasing in $\ell$. By the continuity of $\lm_1(\ell;d_1,a)$, there
exists $\varepsilon>0$ such that $\lm_1(s_\infty;d_1,a-k\varepsilon)<0$.
In view of Theorem \ref{th4.1}, $\dd\lim_{t\to\infty}\|v(t,\cdot)\|_{C([0,s(t)])}=0$.
There exists $\tau\gg 1$ such that $\lm_1(s(\tau);d_1,a-k\varepsilon)<0$ and
$v(t,x)\leq \varepsilon$ for all $t\geq\tau,\,x\in[0,s(\tau)]$. Let $w$ be the
unique solution of
  $$\left\{\begin{array}{ll}
 w_t-d_1w_{xx}=w\big(a(t,x)-k\varepsilon-w\big), \ \ &t>\tau, \
 0<x<s(\tau),\\[1mm]
 B_1[w](t,0)=w(t,s(\tau))=0,&t>\tau,\\[1mm]
  w(\tau,x)=u(\tau,x),&0\le x\le s(\tau).
  \end{array}\right.$$
Then $u\geq w$ in $[\tau,\infty)\times[0,s(\tau)]$ by the comparison principle. Note
that $\lm_1(s(\tau);d_1,a-k\varepsilon)<0$, it follows from Theorem 28.1 of
\cite{Hess} that $w(t+nT,x)\to Z(t,x)$ as $n\to\infty$ uniformly on
$[0,T]\times[0,s(\tau)]$, where $Z(t,x)$ is the unique positive solution of the
following $T$-periodic boundary value problem
\bess\left\{\begin{array}{ll}
 Z_t-d_1Z_{xx}=Z\big(a(t,x)-k\varepsilon-Z\big),\ \ &0\le t\le
 T, \ \ 0<x<s(\tau),\\[1mm]
  B_1[Z](t,0)=Z(t,s(\tau))=0, &0\le t\le T,\\[1mm]
 Z(0,x)=Z(T,x), &0\le x\le s(\tau).\end{array}\right.
 \eess
Since $u\geq w$ in $[\tau,\infty)\times[0,s(\tau)]$, we immediately obtain
 \[\liminf_{n\to\infty}u(t+nT,x)\ge Z(t,x)\ff (t,x)\in
 [0,T]\times[0,s(\tau)].\]
This is a contradiction with the first formula of (\ref{4.1}). The proof is
complete. \ \ \ \fbox{}

In the following, with the parameter $s_0$ satisfying $s_0<s^*$ and $(u_0,v_0)$
being fixed, let us discuss the effect of the coefficient $\mu$ on the spreading and
vanishing. We first give a lemma.
 \begin{lem}\lbl{l5.3}\, Let $d,C>0$ be fixed constants and the boundary operator $B$
be as above. For any given constants $s_0, \Lambda>0$, and any function $\bar u_0\in
C^2([0,s_0])$ satisfying $B[\bar u_0](0)=\bar u_0(s_0)=0$ and $\bar u_0>0$ in
$(0,s_0)$, there exists $\mu^0>0$ such that when $\mu\geq\mu^0$ and $(\bar u, \bar
s)$ satisfies
 \bess
 \left\{\begin{array}{ll}
   \bar u_t-d\bar u_{xx}\geq -C \bar u, \ &t>0, \ 0<x< \bar s(t),\\[1mm]
  B[\bar u](t,0)=0=\bar u(t, \bar s(t)),\ &t\geq 0,\\[1mm]
 \bar s'(t)=-\mu \bar u_x(t, \bar s(t)), \ &t\geq 0,\\[1mm]
  \bar s(0)=s_0, \ \bar u(0,x)=\bar u_0(x),\ \  &0\leq x\leq s_0,
 \end{array}\right.
 \eess
we must have $\dd\liminf_{t\to\infty}\bar s(t)>\Lambda$.
\end{lem}

The proof of Lemma \ref{l5.3} is essentially similar to that of Lemma 3.2 in
\cite{WZjdde} and is hence omitted.

Recalling the estimate (\ref{2.2}), as a consequence of Lemmas \ref{l5.2} and
\ref{l5.3}, we have
  \begin{cor}\lbl{c5.1}\, Assume that ${\cal E}\not=\emptyset$ and set $s^*=\min{\cal E}$. If
$s_0<s^*$, then there exists $\mu^0>0$ depending on $(u_0,v_0,s_0)$ such that
$s_\infty=\infty$ if $\mu\geq\mu^0$.
\end{cor}
 \begin{lem}\label{l5.4}\,Assume that ${\cal E}\not=\emptyset$ and set $s^*=\min{\cal E}$. If
$s_0<s^*$, then there exists  $\mu_0>0$, such that $s_\infty\le s^*$ for all
$\mu\leq\mu_0$.
\end{lem}

{\bf Proof}.\,This proof is similar to that of Lemma 5.2 in \cite{Wperiodic15}. Here we
give the sketch for completeness and readers' convenience. Since $s_0<s^*=\min{\cal E}$,
we have that $\lm_1(s_0;d_1,a)>0$ and $\gamma_1(s_0;d_2,b)>0$.

Let $\phi(t,x)$ and $\psi(t,x)$ be, respectively, the positive eigenfunctions corresponding to
$\lambda_1:=\lm_1(s_0;d_1,a)$ and $\gamma_1:=\gamma_1(s_0;d_2,b)$ of (\ref{3.9}) and
(\ref{3.10}) with $\ell=s_0$. The following conclusions are obvious:

(i)\, $\phi_x(t,s_0)<0$, $\psi_x(t,s_0)<0$ in $[0,T]$;

(ii)\, $\phi(t,0)>0$ in $[0,T]$ when $\beta_1>0$, $\psi(t,0)>0$ in $[0,T]$ when
$\beta_2>0$;

(iii)\, $\phi_x(t,0)>0$ in $[0,T]$ when $\beta_1=0$, and $\psi_x(t,0)>0$ in $[0,T]$
when $\beta_2=0$.

In view of the above facts (i)-(iii) and the regularity of $\phi$ and $\psi$, we
know that there exists a constant $C>0$ such that
 \bes
 x\phi_x(t,x)\le C\phi(t,x),\ \ x\psi_x(t,x)\le C\psi(t,x), \ \ \forall \
 (t,x)\in[0,T]\times[0, s_0].
 \lbl{5.1}\ees
Let $0<\delta,\,\sigma<1$ and $\Lambda>0$ be constants, which will be determined
later. Set
 \bess
 &\dd g(t)=1+2\delta-\delta {\rm e}^{-\sigma t}, \ \ \xi(t)=\int_0^t
 g^{-2}(\tau){\rm d}\tau, \ \ t\geq 0,& \\[1mm]
 &w(t,x)=\Lambda{\rm e}^{-\sigma t}\phi\left(\xi(t),\,y\right), \ \
 z(t,x)=\Lambda{\rm e}^{-\sigma t}\psi\left(\xi(t),\,y\right), \ \ y=\dd\frac x{g(t)},
 \ \ 0\leq x\leq s_0g(t).&
 \eess

Firstly, for any given $0<\varepsilon\ll 1$, since $a$ and $b$ are uniformly
continuous in $[0,T]\times[0,3s_0]$ and $T$-periodic in $t$, there
exists $0<\delta_0(\ep)\ll 1$ such that, for all $0<\delta\le\delta_0(\ep)$ and
$0<\sigma<1$,
 \bes
 \left|a\big(\xi(t),\,y(t,x)\big)-g^2(t)a(t,x)\right|\le\varepsilon,\ \
 \left|b\big(\xi(t),\,y(t,x)\big)-g^2(t)b(t,x)\right|\le\varepsilon
  \lbl{5.2}\ees
for $t\ge 0$ and $0\leq x\leq s_0g(t)$. Remembering $\lambda_1,\,\gamma_1>0$,
in view of (\ref{5.1}) and (\ref{5.2}), the direct calculation yields that
 \bes
 w_t-d_1w_{xx}-w\big(a(t,x)-w\big)&>&v(-\sigma-\varepsilon-\sigma C+\lambda_1/4)
 >0,\lbl{5.3} \\[1mm]
 z_t-d_2z_{xx}-z\big(b(t,x)-z\big)&>&z(-\sigma-\varepsilon-\sigma C+\gamma_1/4)
 >0
 \lbl{5.4}\ees
for all $t>0$ and $0<x<s_0g(t)$ provided $0<\sigma,\varepsilon\ll 1$. Carefully
analysis gives
  \bes
 B_1[w](t,0)\ge 0, \ \ B_2[z](t,0)\ge 0, \ \ w(t,s_0g(t))=0, \ \ z(t,s_0g(t))=0, \ \
 \forall \  t>0.
 \lbl{5.5}\ees
Fix $0<\sigma,\varepsilon\ll 1$ and $0<\delta\le\delta_0(\ep)$. Based on the
regularities of $u_0(x),\,v_0(x),\,\phi(0,x)$ and $\psi(0,x)$, one can choose
$\Lambda\gg1$ such that
 \bes
 u_0(x)\leq \Lambda\phi\left(0,\,\frac x{1+\delta}\right)=w(0,x), \ \ v_0(x)\leq
 \Lambda\psi\left(0,\,\frac x{1+\delta}\right)=z(0,x), \ \ \forall \ 0\le x\le s_0.
\lbl{5.6}\ees
According to $s_0g'(t)=s_0\sigma\delta{\rm e}^{-\sigma t}$ and
 \[w_x(t,s_0g(t))=\frac 1{g(t)}\Lambda{\rm e}^{-\sigma t}\phi_y(\xi(t),s_0),\ \
  z_x(t,s_0g(t))=\frac 1{g(t)}\Lambda{\rm e}^{-\sigma t}\psi_y(\xi(t),s_0),\]
there exists $\mu_0>0$ such that
 \bes
 s_0g'(t)\geq-\mu\big(w_x(t,s_0g(t))+\rho z_x(t,s_0g(t))\big), \ \ \forall \
 0<\mu\le \mu_0, \ t\ge 0.
\lbl{5.7}\ees

Remembering (\ref{5.3})-(\ref{5.7}), we can apply the comparison principle (Lemma \ref{l5.1}) to $(u,v, s(t))$ and
$(w,z,s_0g(t))$, and then derive that
 \[s(t)\leq s_0g(t), \ \ u(t,x)\leq w(t,x), \ \ v(t,x)\leq z(t,x), \ \ \forall \
 t\geq 0, \ 0\leq x\leq s(t).\]
Hence $s_\infty\leq s_0g(\infty)=s_0(1+2\delta)$ for all $0<\mu\leq\mu_0$. The proof
is complete.  \ \ \ \fbox{}

Now, let us give the criteria for spreading and vanishing of the problem
(\ref{1.2}).
  \begin{theo}\lbl{th5.1} Assume that ${\cal E}\not=\emptyset$ and set $s^*=\min{\cal E}$.

{\rm(i)} If $s_0\ge s^*$, then $s_\infty=\infty$ for all $\mu>0$;

{\rm(ii)}\, If $s_0<s^*$, then  there exist $\mu^*\geq\mu_*>0$, depending on
$(u_0,v_0,s_0)$, such that $s_\infty\leq s^*$ if $\mu\leq\mu_*$, and
$s_\infty=\infty$ if $\mu>\mu^*$.
\end{theo}

{\bf Proof}. Remember Lemmas \ref{l5.2} and \ref{l5.4} and Corollary \ref{c5.1},
the proof is similar to that of Theorem 5.2 in \cite{Wjde14}.
We omit the details.\ \ \ \fbox{}

We have mentioned in the above that if one of the functions $a(t,x)$ and $b(t,x)$ satisfies the condition {\bf(A)}, then ${\cal E}\not=\emptyset$.

\section{The problem (\ref{1.3})}\lbl{s5}
\setcounter{equation}{0}

In this section, we briefly discuss the problem (\ref{1.3}). So, $(u,v,s)$ means the unique solution of {\rm(\ref{1.3})} throughout this section.
  \begin{theo}\lbl{th6.2}{\rm(}Vanishing{\rm)}\, Assume that the condition {\bf(H1)}
holds. If $s_\infty<\infty$, then
 \bes
 &\dd\lim_{t\to\infty}\|u(t,\cdot)\|_{C([0,s(t)])}=0,\lbl{6.1}&\\[1mm]
 &\dd\lim_{n\to\infty}v(t+nT,x)=V(t,x) \ \ \mbox{uniformly\ in } \
 [0,T]\times[0,L]&\lbl{6.2}
 \ees
for any $L>0$, where $V$ is the unique positive solution of the following $T$-periodic boundary value problem
  \bes\left\{\begin{array}{ll}
 V_t-d_2V_{xx}=V\big(b(t,x)-V\big), \ \ &0\le t\le T, \
 0<x<\yy,\\[1mm]
 B_2[V](t,0)=0,  &0\le t\le T,\\[1mm]
 V(0,x)=V(T,x), &0\le x<\yy.
\end{array}\right.\lbl{6.2a}\ees
This shows that if the invasive species $u$ cannot spread successfully, it will extinct in
the long run.
\end{theo}

{\bf Proof}. The proof of (\ref{6.1}) is the same as that of (\ref{4.1}). The
existence and uniqueness of $V(t,x)$ is guaranteed by Proposition \ref{p3.3}.
Thanks to (\ref{6.1}), similarly to the proof of Theorem \ref{th4.2}, we can show
that
$\dd\liminf_{n\to\infty}v(t+nT,x)\geq V(t,x)$ and
$\dd\limsup_{n\to\infty}v(t+nT,x)\leq V(t,x)$ uniformly in $[0, T]\times[0,L]$.
This finishes the proof.\ \ \ \fbox{}

\vskip 2pt When $s_\infty=\infty$, Theorem \ref{th4.2} (Spreading) remains hold.

In the following three lemmas, we assume that the function $a(t,x)$ satisfies condition
{\bf(A)} and take $s_*>0$ such that $\lm_1(s_*;d_1,a)=0$. Obviously, $s_*$ exists
uniquely.
  \begin{lem}\lbl{l6.1}\, If $s_\infty<\infty$, then $s_\infty\leq s_*$. Hence,
$s_0\ge s_*$ implies $s_\infty=\infty$ for all $\mu>0$.
\end{lem}

The proof of Lemma \ref{l6.1} is essentially same as that of Lemma \ref{l5.2} and
is hence omitted here.
  \begin{lem}\lbl{l6.2}\,If $s_0<s_*$, then there exist $0<\mu_0<\mu^0$, such that
$s_\infty=\infty$ for  $\mu\geq\mu^0$, and $s_\infty\leq s_*$ for $\mu\le\mu_0$.
\end{lem}

{\bf Proof}. Recalling the estimates obtained in Theorem \ref{th6.1}, and applying Lemmas \ref{l5.3} and \ref{l6.1}, we can derive that there exists $\mu^0>0$, such that
$s_\infty=\infty$ for $\mu\geq\mu^0$. On the other hand, since $u$ satisfies
 \bess\left\{\begin{array}{lll}
 u_t-d_1u_{xx}<u\big(a(t,x)-u\big), &t>0, \ 0<x<s(t),\\[1mm]
  B_1[u]=0,\ \ &t>0, \ =0,\\[1mm]
 u=0, \ s'(t)=-\mu u_x,\ \ &t>0, \ x=s(t), \\[1mm]
 s(0)=s_0, \ \ u(0,x)=u_0(x), \ & 0\leq x\leq s_0,
 \end{array}\right.\eess
by the known results for the logistic equation (see Lemma 5.2 of \cite{Wperiodic15}) and
comparison principle, we can show that there exists $\mu_0>0$, such that
$s_\infty\leq s_*$ for $\mu\le\mu_0$. The proof is complete. \ \ \ \fbox{}

In the same way as the proof of Lemma 2.6 in \cite{DL2}, it can be proved that $(u,
v, s)$ is monotone increasing in $\mu$. Similarly to the proof of Lemma 4.9 in
\cite{DL2}, we have the following lemma.
  \begin{lem}\lbl{l6.3}\, If $s_0<s_*$, then there exist $\mu^*>0$, such that
$s_\infty=\infty$ for  $\mu>\mu^*$, while $s_\infty\leq s_*$ for $\mu\le\mu_*$.
\end{lem}

It is worth mentioning that the assumption {\bf(H1)} implies condition {\bf(A)}.
Summarizing the above conclusions we obtain the following spreading-vanishing
dichotomy and sharp criteria for spreading and vanishing.
 \begin{theo}\lbl{th6.3}\, Under the condition {\bf(H1)}, we have the following alternative conclusion:

Either\vspace{-2mm}\begin{itemize}
\item[{\rm(i)}]spreading: $s_\infty=\infty$, both $(\ref{4.2})$ and $(\ref{4.3})$ hold,
  \vspace{-2mm}\end{itemize}
or
\vspace{-2mm}\begin{itemize}
\item[{\rm(ii)}]
vanishing: $s_\infty\le s_*$, and
$\dd\lim_{t\to\infty}\|u(t,\cdot)\|_{C([0,s(t)])}=0$,
$\dd\lim_{n\to\infty}v(t+nT,x)=V(t,x)$ uniformly in $[0,T]\times[0,L]$ for any given
$L>0$, where $V$ is the unique positive $T$-periodic solution of $(\ref{6.2a})$.
\end{itemize}\vspace{-1mm}

Moreover,
\vspace{-2mm}\begin{itemize}
\item[{\rm(iii)}]
If $s_0\ge s_*$, then $s_\infty=\infty$ for all $\mu>0$.\vspace{-1mm}

\item[{\rm(iv)}] If $s_0<s_*$, then there exist $\mu^*>0$, such that $s_\infty=\infty$
for $\mu>\mu^*$, while $s_\infty\leq s_*$ for $\mu\leq\mu^*$.\end{itemize}\vspace{-2mm}
\end{theo}

Finally, we estimate the asymptotic spreading speed of the free boundary $s(t)$ when
spreading occurs. To this aim, let us first state a known result, which plays an
important role in the study of asymptotic spreading speed. For a $T$-periodic function $f(t)$, we define
   \[\od f=\frac 1T\int_0^Tf(t){\rm d}t.\]
  \begin{prop}\lbl{p6.1}{\rm(\cite[Section 2]{DGP})}\, Let $d>0$ and $0<\nu<1$ be the given constants. Assume that $F, \varphi\in C^\nu([0,T])$ are $T$-periodic functions, $\varphi$ is  positive and $F$ is nonnegative in $[0,T]$. Then the problem
\bess\left\{\begin{array}{ll}
 w_t-dw_{xx}+F(t)w_x=\varphi(t)w-w^2, \ \ &0\le t\le T, \ 0<x<\yy,\\[1mm]
 w(t,0)=0, &0\le t\le T,\\[1mm]
w(0,x)=w(T,x),&0\le x<\yy
  \end{array}\right.\eess
has a positive $T$-periodic solution $w^F(t,x)\in C^{1,2}([0,T]\times[0,\infty))$ if and only if  $\od F<2\sqrt{d\od\vp}$, and such a solution is unique when it exists.  Moreover, the following hold: 

{\rm(i)} $w^F_x(t,x)>0$ and $w^F(t,x)\to z(t)$ uniformly in $[0,T]$ as $x\to\infty$, where
$z(t)$ is the unique positive periodic solution of the problem
 \[z'=\varphi(t)z-z^2, \ \ 0\le t\le T; \ \ \ z(0)=z(T);\]

{\rm(ii)} For any given nonnegative $T$-periodic function $G\in
C^\nu([0,T])$ satisfying  $\od G<2\sqrt{d\od\vp}$, the assumption $G\le,\,\not\equiv F$ implies
$w^G_x(t,0)>w^{F}_x(t,0)$, $w^G(t,x)>w^{F}(t,x)$ for $0\le t\le T$ and $x>0$;

{\rm(iii)} For each $\mu>0$, there exists a unique positive $T$-periodic function $F_0(t)=F_0(d,\mu,\varphi)(t)\in C^\nu([0,T])$ such that $\mu
w^{F_0}_x(t,0)=F_0(t)$ in $[0,T]$,  and $0<\od F_0<2\sqrt{\od \varphi d}$.
 \end{prop}
  \begin{theo}\lbl{th6.4} Under the condition {\bf(H1)} with $r=0$ we further assume that $(\ref{3.11})$ hold. When
the spreading occurs, i.e., $s_\infty=\infty$, we have
 \bess
 \limsup_{t\to\infty}\frac{s(t)}t\le
 \frac 1T\int_0^T F_0(d_1,\mu,a^\infty-k z_1)(t){\rm d}t,\ \
 \liminf_{t\to\infty}\frac{s(t)}t\ge \frac 1T\int_0^T F_0(d_1,\mu,a_\infty-k
 z_2)(t){\rm d}t,
 \eess
where $z_1(t)$ and $ z_2(t)$ are the unique positive solutions of $(\ref{3.14})$ and
$(\ref{3.15})$, respectively.
 \end{theo}

{\bf Proof}.\, This proof is similar to that of Theorem 4.4 in \cite{DGP} with some
modifications. Here we give the sketch for completeness and readers' convenience.

First of all, in view of (\ref{3.14}) and (\ref{3.15}), it is easy to see that
 \[z_1(t)\le\max_{[0,T]}b_\yy(t)\le\bar b^\yy, \ \ z_2(t)\le\max_{[0,T]}b^\yy(t)=\bar b^\yy, \ \ \forall \ t\in[0,T].\]
Therefore,
 \[a^\infty(t)-k z_1(t)\ge \ud a_\yy-k\bar b^\yy>0, \ \ a_\infty(t)-k z_2(t)\ge \ud a_\yy-k\bar b^\yy>0, \ \ \forall \ t\in[0,T]\]
by the condition  (\ref{3.11}).

{\bf Step 1}.\, Let  $(U^*,V_*)$ and $(U_*,V^*)$ be positive
solutions of $(\ref{3.1})$ obtained in Theorem \ref{th3.1}. Apply the last two
conclusions of (\ref{3.13}), it yields that
 \bess
 z_1(t)\le\liminf_{x\to\infty} V_*(t,x),\ \ \
 \limsup_{x\to\infty} V^*(t,x)\le z_2(t)
 \eess
uniformly in $[0,T]$. By the same way as that of step 1 in the proof of
\cite[Theorem 4.4]{DGP}, we can show that
 \bes
 \limsup_{x\to\infty}U^*(t,x)\le\od \psi(t), \ \ \
 \liminf_{x\to\infty}U_*(t,x)\ge\ud \psi(t) \ \ \ \mbox{uniformly in}\ \
 [0,T],
 \lbl{6.4}\ees
where $\od\psi(t)$ and $\ud \psi(t)$ are, respectively, the unique positive
$T$-periodic solutions of
 \[\od\psi'(t)=\od \psi\big(a^\infty(t)-k z_1(t)-\od \psi\big), \ \ \ \od
 \psi(0)=\od \psi(T)\]
and
 \[\ud\psi'(t)=\ud \psi\big(a_\infty(t)-k z_2(t)-\ud \psi\big), \ \ \ \ud
 \psi(0)=\ud \psi(T).\]

{\bf Step 2}.\, For the given $0<\varepsilon\ll 1$, by (\ref{6.4}), there exists
$\ell^*=\ell^*(\varepsilon)\gg 1$ such that
 \[U^*(t,x)\le\od\psi_\varepsilon(t),\ \ U_*(t,x)\ge\ud\psi_\varepsilon(t)
 \ \ \ \mbox{in}\ \
 [0,T]\times[\ell^*,\infty),\]
where $\od \psi_\varepsilon(t)$ and $\ud \psi_\varepsilon(t)$ are, respectively, the
unique positive $T$-periodic solutions of
 \[\od\psi'_\varepsilon(t)=\od\psi_\varepsilon\big(a^\infty(t)-k
 z_1(t)+\varepsilon-\od\psi_\varepsilon\big), \ \ \
 \od\psi_\varepsilon(0)=\od\psi_\varepsilon(T)\]
and
 \[\ud\psi'_\varepsilon(t)=\ud\psi_\varepsilon\big(a_\infty(t)-k
 z_2(t)-\varepsilon-\ud\psi_\varepsilon\big), \ \ \
 \ud\psi_\varepsilon(0)=\ud\psi_\varepsilon(T).\]
Because Theorem \ref{th4.2} holds for the problem (\ref{1.3}), in view of
(\ref{4.2}), (\ref{4.3}) and $s_\infty=\infty$, we have that there exists a positive
integer $n=n(\ell^*)$ such that $s(nT)>3\ell^*$ and
 \[\ud\psi_{2\varepsilon}(t)\le u(t+nT,2\ell^*)\le\od\psi_{2\varepsilon}(t)\ff t\ge
 0.\]
Follow the arguments of steps 2 and 3 in the proof of \cite[Theorem 4.4]{DGP} step
by step, we can obtain the desired results. The details are omitted here. \ \ \
\fbox{}
  \begin{remark}\lbl{r6.1}  The main difference between $(\ref{1.2})$ and
$(\ref{1.3})$  is the following. For the problem $(\ref{1.2})$, the criteria we got for
spreading and vanishing are not sharp, due to the lack of monotonicity of solution
in $\mu$. Also, the estimate of the asymptotic spreading speed of free boundary has
not been obtained for the problem $(\ref{1.2})$.
\end{remark}


\begin{thebibliography}{99}
\bibliographystyle{siam}
\setlength{\baselineskip}{15pt}

\vspace{-1.5mm}\bibitem{ChenA} X.F. Chen \& A. Friedman, {\it A free boundary
problem arising in a model of wound healing}, SIAM J. Math. Anal., {\bf
32}(4)(2000), 778-800.

\vspace{-1.5mm}\bibitem{CC2} R.S. Cantrell \& C. Cosner, {\it Spatial Ecology via
Reaction-Diffusion Equations}, Wiley Series in Mathematical and Computational
Biology, John Wiley \& Sons Ltd, 2003.

\vspace{-1.5mm}\bibitem{CLL} X.F. Chen, K.-Y. Lam \& Y. Lou, {\it Dynamics of a
reaction-diffusion-advection model for two competing species}, Discrete Contin. Dyn.
Syst., Ser. A, {\bf 32}(2012), 3841-3859.

\vspace{-1.5mm}\bibitem{DHMP} J. Dockery, V. Hutson, K. Mischaikow \& M.
Pernarowski, {\it The evolution of slow dispersal rates: a reaction diffusion
model}, J. Math. Biol.,  {\bf 37}(1998), 61-83.

\vspace{-1.5mm}\bibitem{DG} Y.H. Du \& Z. M. Guo, {\it Spreading-vanishing dichotomy
in
the diffusive logistic model with a free boundary}, II, J. Differential Equations,
{\bf 250}(2011), 4336-4366.

\vspace{-1.5mm}\bibitem{DG1} Y.H. Du \& Z.M. Guo, {\it The Stefan problem for the
Fisher-KPP equation}, J. Differential Equations, {\bf 253}(3)(2012), 996-1035.

\vspace{-1.5mm}\bibitem{DGP} Y.H. Du, Z.M. Guo \& R. Peng, {\it A diffusive logistic
model with a free boundary in time-periodic environment}, J. Funct. Anal., {\bf
265}(2013), 2089-2142.

\vspace{-1.5mm}\bibitem{DLiang} Y.H. Du \& X. Liang, {\it Pulsating semi-waves in
periodic media and spreading speed determined by a free boundary model}, Ann. Inst.
Henri Poincare Anal. Non Lineaire (2013),
http://dx.doi.org/10.1016/j.anihpc.2013.11.004.

\vspace{-1.5mm}\bibitem{DLin} Y.H. Du \& Z.G. Lin, {\it Spreading-vanishing
dichotomy in the diffusive logistic model with a free boundary}, SIAM J. Math.
Anal., {\bf 42}(2010), 377-405.

\vspace{-1.5mm}\bibitem{DL2} Y.H. Du \& Z.G. Lin, {\it The diffusive competition
model with a free boundary: Invasion of a superior or inferior competitor}, Discrete Cont. Dyn. Syst.-B, {\bf 19}(10)(2014), 3105-3132.

\vspace{-1.5mm}\bibitem{DLou} Y.H. Du \& B.D. Lou, {\it Spreading and vanishing in
nonlinear diffusion problems with free boundaries}, J. Eur. Math. Soc., to appear
(arXiv1301.5373)

\vspace{-1.5mm}\bibitem{DMZ} Y.H. Du, H. Matsuzawa \& M.L. Zhou, {\it Sharp estimate
of the spreading speed determined by nonlinear free boundary problems}, SIAM J.
Math. Anal., {\bf 46}(1)(2014), 375-396.

 %\vspace{-1.5mm}\bibitem{DP} Y.H. Du \& R. Peng, {\it The Periodic
 %logistic equation
  %with spatial and temporal degeneracies}, Trans. Amer. Math. Soc., {\bf
 %364}(11)(2012), %6039-6070.

\vspace{-1.5mm}\bibitem{GL} J.S. Guo \& X. Liang,{\it The minimal speed of traveling
fronts for the Lotka-Volterra competition system}, J.
Dyn. Diff. Equat., {\bf 23}(2011), 353-363.

\vspace{-1.5mm}\bibitem{GW} J.S. Guo \& C.H. Wu, {\it On a free boundary problem for
a two-species weak competition system}, J. Dyn. Diff. Equat., {\bf 24}(2012),
873-895.

%\vspace{-1.5mm}\bibitem{HC} A. Hastings et al., {\it The spatial spread of
%invasions: %new developments in theory and evidence}, Ecol. Lett., {\bf 8}(2005),
%91-101.

\vspace{-1.5mm}\bibitem{Hess} P. Hess, {\it Periodic-Parabolic Boundary Value
Problems and Positivity}, Pitman Res. Notes Math., vol. 247, Longman
Sci. Tech., Harlow, 1991.

\vspace{-1.5mm}\bibitem{HMP}V. Hutson, K. Mischaikow \& P. Pol\'{a}\v{c}ik, {\it
The evolution of dispersal rates in a heterogeneous time-periodic environment}, J.
Math. Bio., {\bf 43}(2001), 501-533.

\vspace{-1.5mm}\bibitem{Ka} Y. Kaneko, {\it Spreading and vanishing behaviors for radially symmetric solutions of free boundary problems for reaction-diffusion
equations}, Nonlinear Anal.: Real World Appl., \textbf{18}(2014), 121-140.

\vspace{-1.5mm}\bibitem{KY} Y. Kaneko \& Y. Yamada, {\it A free boundary problem for
a reaction diffusion equation appearing in ecology}, Advan. Math. Sci. Appl., {\bf
21}(2)(2011), 467-492.

%\vspace{-1.5mm}\bibitem{KL1} C. S. Kolar \& D. M. Lodge, {\it Progress in invasion
%biology: %predicting invaders}, Trends Ecol. Evol., {\bf 16}(2001), 199-204.

%\vspace{-1.5mm}\bibitem{KL2} C.S. Kolar \& D.M. Lodge, {\it Ecological predictions
%and risk %assessment for alien fishes in North American}, Science, {\bf 298}(2002),
%1233-1236.

\vspace{-1.5mm}\bibitem{LSU} O. A. Ladyzenskaja, V. A. Solonnikov \& N. N. Uralceva,
{\em Linear and Quasilinear Equations of Parabolic Type}, Academic Press, New York,
London, 1968.

\vspace{-1.5mm}\bibitem{LL} K.Y. Lam \& Y. Lou, {\it Evolution of conditional
dispersal: evolutionarily stable strategies in spatial models}, J. Math. Bio., {\bf
68}(2014), 851-877.

\vspace{-1.5mm}\bibitem{LLW} M.A. Lewis, B.T. Li \& H.F. Weinberger, {\it Spreading
speed and linear determinacy for two-species competition models}, J. Math. Biol.,
{\bf 45}(2002), 219-233.

%\vspace{-1.5mm}\bibitem{LZ} X. Liang \& X.Q. Zhao, {\it Asymptotic speeds of spread
%and traveling waves for monotone semiflows with
%applications}, Commun. Pure Appl. Math., {\bf 60}(2007), 1-40.

\vspace{-1.5mm}\bibitem{Lie} G. M. Lieberman, {\it Second Order Parabolic
Differential Equations}, World Scientific Publishing Co. Inc., River Edge, NJ, 1996.

\vspace{-1.5mm}\bibitem{LHP} J.L. Lockwood, M.F. Hoopes \& M.P. Marchetti, {\it
Invasion Ecology}, Blackwell Publishing, 2007.

\vspace{-1.5mm}\bibitem{Ni} W.M. Ni, {\it The Mathematics of Diffusion}, CBMS-NSF
Regional Conf. Ser. in Appl. Math. 82, SIAM, Philadelphia, 2011.

\vspace{-1.5mm}\bibitem{PD} R. Peng \& W. Dong, {\it The periodic-parabolic logistic
equation on $\mR^N$}, Discrete Cont. Dyn. Syst. A, {\bf 32}(2)(2012),
619-641.

\vspace{-1.5mm}\bibitem{PZ} R. Peng \& X.Q. Zhao, {\it  The diffusive logistic model
with a free boundary and seasonal succession}, Discrete Cont. Dyn. Syst. A,
{\bf 33}(5)(2013), 2007-2031.

\vspace{-1.5mm}\bibitem{PW} M.H. Protter \& H.F. Weinberger, {\it Maximum
Principles in Differential Equations}, Prentice Hall: Englewood
Cliffs, 1967.

\vspace{-1.5mm}\bibitem{SK} N. Shigesada \& K. Kawasaki, {\it Biological Invasions:
Theory and Practice}, Oxford Ser. Ecol. Evol., Oxford Univ. Press,
Oxford, 1997.

\vspace{-1.5mm}\bibitem{WZ15} J. Wang and L. Zhang, {\it Invasion by an inferior or superior competitor: A diffusive competition model with a free boundary in a heterogeneous environment}, J. Math. Anal. Appl., {\bf 423}(2015), 377-398.

\vspace{-1.5mm}\bibitem{Wjde14} M.X. Wang, {\it On some free boundary problems of the
prey-predator model}, J. Differential Equations, {\bf 256}(10)(2014), 3365-3394.

\vspace{-1.5mm}\bibitem{Wjde15} M.X. Wang, {\it The diffusive logistic equation with
a free boundary and sign-changing coefficient}, J. Differential Equations, {\bf 258}(2015), 1252-1266.

\vspace{-1.5mm}\bibitem{Wcnsns15} M.X. Wang, {\it Spreading and vanishing in the diffusive prey-predator model with a free boundary},  Commun. Nonlinear Sci. Numer. Simulat., {\bf 23}(2015), 311-327.

\vspace{-1.5mm}\bibitem{Wperiodic15} M.X. Wang, {\it A diffusive logistic equation with
a free boundary and sign-changing coefficient in time-periodic environment}.
arXiv:1504.03958 [math.AP], 2015.

\vspace{-1.5mm}\bibitem{WZhang} M.X. Wang and Y. Zhang, {\it Two kinds of free boundary problems for the diffusive prey-predator model}, Nonlinear Anal.: Real World Appl., {\bf 24}(2015), 73-82.

\vspace{-1.5mm}\bibitem{WZ} M.X. Wang \& J.F. Zhao, {\it A free boundary problem for a predator-prey model with double free boundaries}. arXiv:1312.7751 [math.DS], 2013.

\vspace{-1.5mm}\bibitem{WZjdde} M.X. Wang \& J.F. Zhao, {\it Free boundary problems
for a Lotka-Volterra competition system}, J. Dyn. Diff. Equat., {\bf 26}(3)(2014), 655-672 (DOI: 10.1007/s10884-014-9363-4)

\vspace{-1.5mm}\bibitem{YLWW} Q.X. Ye, Z.Y. Li, M.X. Wang \& Y.P. Wu, {\it
Introduction to Reaction-Diffusion Equations} (in chinese), Science Press, Beijing,
2011.

\vspace{-1.5mm}\bibitem{ZW} J.F. Zhao \& M.X. Wang, {\it A free boundary problem
of a predator-prey model with higher dimension and heterogeneous environment},
Nonlinear Anal.: Real World Appl.,  {\bf 16}(2014), 250-263.

\vspace{-1.5mm}\bibitem{Zhaox} X.Q. Zhao, {\it Spatial dynamics of some evolution
systems in biology}, in Y.H. Du, H. Ishii \& W. Y. Lin (eds.): Recent Progress on Reaction-Diffusion Systems and Viscosity Solutions, World Scientific, Singapore, 2009.

\vspace{-1.5mm}\bibitem{ZX} P. Zhou \& D.M. Xiao, {\it The diffusive logistic model
with a free boundary in heterogeneous environment}, J. Differential Equations, {\bf
256}(2014),  1927-1954.

\end{thebibliography}
\end{document}